\documentclass[12pt]{article}
\usepackage{amsmath, amssymb}
\begin{document}
\centerline{\bf ABUNDANCE CONJECTURE}

\bigbreak\centerline{\it Dedicated to Shing-Tung Yau}

\bigbreak\centerline{Yum-Tong Siu\ %
\footnote{Partially supported by Grant 0500964 of the National Science
Foundation.  Written for the Festschrift of the sixtieth birthday of Shing-Tung Yau} }

\bigbreak

\bigbreak\noindent{\bf\S0. Introduction.}  For a compact complex algebraic manifold $X$ we define its Kodaira dimension $$\kappa_{\rm kod}(X)=\limsup_{m\to\infty}\frac{\log\dim_{\mathbb
C}\Gamma\left(X,mK_X\right)}{\log m}$$ and its numerical Kodaira dimension $$\kappa_{\rm num}(X)=\sup_{k\geq 1}\left[\limsup_{m\to\infty}\frac{\log\dim_{\mathbb
C}\Gamma\left(X,mK_X+kA\right)}{\log m}\right],$$
where $A$ is any ample line bundle on $X$.  The definition of $\kappa_{\rm num}(X)$ is independent of the choice of the ample line bundle $A$ on $X$.
We use the analytic viewpoint and formulate the abundance conjecture as $\kappa_{\rm kod}(X)=\kappa_{\rm num}(X)$.  For terminology and formulations in the context of algebraic geometry see [Kawamata1985, p.569] and [Nakayama2004] and the remarks at the end of Chapter 7 of [Koll\'ar-Mori1998].  In this note we present a sketch of the proof of the following theorem which affirms the abundance conjecture.  Complete details are being written up and will be available later elsewhere.

\bigbreak\noindent(0.1) {\it Theorem.} For a compact complex algebraic manifold $X$ the Kodaira dimension $\kappa_{\rm kod}(X)$ of $X$ equals the numerical Kodaira dimension $\kappa_{\rm num}(X)$ of $X$.

\bigbreak In this sketch some parts are given in details while for some others only the key techniques and steps are explained.  Even though in the parts where only a sketch is given, enough of the essential arguments is presented for the reader to understand how the proof works.  The emphasis of our presentation is to highlight what the important arguments are and why they work and to minimize the amount of prerequisite background material.

\bigbreak\noindent(0.2) {\it Three Parts of the Proof.} The proof is divided into the following three parts.  The first part is the special case of zero numerical Kodaira dimension.  The second part is the general case of positive numerical Kodaira dimension under the additional assumption that the numerically trivial foliation for the canonical line bundle coincides with the numerically trivial fibration for the canonical line bundle.  In other words, the assumption is that the leaves of the numerically trivial foliation for the canonical line bundle are knonwn to be compact.  The third part verifies that the numerically trivial foliation for the canonical line bundle always coincides with the numerically trivial fibration for the canonical line bundle.

\medbreak In all three parts the arguments work only for the canonical line bundle and do not work for general holomorphic line bundles for the reasons which can be very briefly described as follows. For the first part the Kodaira-Serre duality is used to identify the space of all flatly twisted canonical sections with the dual of a cohomology group with coefficients in a flat bundle.  In the second part the variation of Hodge structure is used.  In the third part the adjunction for foliation is used.

\medbreak Besides the usual analytic techniques of $L^2$ estimates of $\bar\partial$ and the standard techniques of algebraic geometry, for the first part we use the technique developed first by Gelfond and Schneider [Gelfond1934, Schneider1934] for their solution of Hilbert's seventh problem.  The technique of Gelfond and Schneider comes in as an ingredient of the method of Simpson [Simpson1993] to replace the flat line bundle in a nontrivial flatly twisted canonical section by one which is a torsion flat line bundle.  For the second part we use the technique of the semipositivity, and the condition for positivity, of the zeroth direct image of the relative pluricanonical line bundle.  For the third part we use the technique of the First Main Theorem of Nevanlinna theory which is related to the technique of Gelfond and Schneider in the first part.

\bigbreak\noindent(0.3) {\it Algebraic Properties of Solutions of Algebraically Defined Differential Equations.}  It is interesting to point out that the special techniques used for the first and second parts can be considered as originating from the same theme of algebraic properties of solutions of algebraically defined differential equations.  Moreover, the technique of the third part is related to that of the first part.  By an algebraically defined differential equation is meant a differential equation whose coefficients are algebraic functions (which may also be required to be defined over an algebraic number field).   The technique of Gelfond-Schneider [Gelfond1934, Schneider1934] concerns the problem of algebraic values of solutions of algebraically defined differential equations.

\medbreak The technique of semipositivity of the zeroth direct image of the relative canonical line bundle can be traced to the 1873 paper of Schwarz [Schwarz1873] where he studied the problem of Gauss about the condition for the existence of algebraic solutions for the Gauss hypergeometric differential equation [Gauss1812].  After the work of Newton on laws of motion and gravitation the question of algebraic solutions of singular differential equations such as the Gauss hypergeometric differential equation is significant because the motions of heavenly bodies satisfy second-order differential equations with singularity coming from mass and algebraic solutions are considered in order to study periodic solutions.

\medbreak Schwarz [Schwarz1873] studied the problem by using Euler's integral representation of the hypergeometric function by beta functions [Euler1778] and interpreting it as the period of an elliptic curve so that the Gauss hypergeometric differential equation is the differential equation from the variation of the period (which in modern terminology means the variation of Hodge structure) with the differentiation used in the differential equation defined by the Gauss-Manin connection.  Actually in the situation considered by Schwarz it is the equivariant variation of Hodge structure with cyclic group action rather than the usual variation of Hodge structure.

\medbreak The zeroth direct image of the relative canonical line bundle is a holomorphic subbundle of the flat bundle whose fibers are the Betti groups endowed with the natural metric defined by the integration of the exterior product of the representative differential form with its complex conjugate.  The curvature of a holomorphic subbundle of a flat vector bundle is seminegative if the metric of the ambient flat vector bundle is positive definite, otherwise the curvature may be semipositive if the second fundamental form has value in the part of the ambient flat vector where the metric is negative definite.

\medbreak In the 1973 paper of Schwarz he did not consider the curvature and the second fundamental form.  Instead he considered the monodromy which he showed in his case to be in the biholomorphism group of an open $1$-disk ({\it i.e.,} elements of $U(1,1)$), which means that in his case the second fundamental form has value in the part of the ambient flat vector bundle where the metric is negative.  In Picard's presentation [Picard1885] of the method of Schwarz and its generalization to the case of complex dimension two the r\^ole of the invariant circle [Picard1885, p.365, Equation(C)] and the r\^ole of $U(2,1)$ [Picard1885, p.378] are more clearly spelled out.

\medbreak With the benefit of hindsight of more than a century later one can say that the seed is already in the 1973 paper of Schwarz for the statement that the curvature of the zeroth direct image of the relative canonical line bundle is semipositive and is given by the second fundamental form in the flat vector bundle whose fibers are the Betti groups.  The papers of Deligne-Mostow [Deligne-Mostow1986] and Mostow [Mostow1987] gave a historic account of the developments from the work of Euler, Schwarz, Picard, Levavasseur [LeVavasseur1893] and their own work from the viewpoint of the monodromy of hypergeometric functions and nonlattice integral monodromy.

\bigbreak\noindent(0.4) {\it Notations.}  For complex manifold $Y$ we denote by $T_Y$ its tangent bundle (of type (1,0)), denote by $T_Y^*$ the dual of $T_Y$, denote by $K_Y$ its canonical line bundle, and denote by $\Omega^p_Y$ the bundle of $(p,0)$-forms on $Y$.  For a holomorphic family $\pi:Y\to S$ of complex manifolds over a complex manifold, we denote by $\Omega^p_{Y|S}$ the relative bundle of $(p,0)$-forms which is the vector bundle on $Y$ whose restriction to a fiber $Y_s$ is equal to the bundle $\Omega_{Y_s}^p$ of $(p,0)$-forms on the fiber $Y_s$ and we denote by $K_{Y|S}$ the relative canonical line bundle which is the line bundle on $Y$ whose restriction to a fiber $Y_s$ is equal to the canonical line bundle $K_{Y_s}$ of the fiber $Y_s$.  The topological boundary of a set $E$ is denoted by $\partial E$ and its topological closure is denoted by $\bar E$.  If there is no confusion, for a divisor $D$ in a complex manifold we also use the same notation $D$ to denote the holomorphic line bundle on $Y$ associated to the divisor $D$ and we use $s_D$ to denote the holomorphic section of the holomorphic line bundle $D$ whose divisor is the divisor $D$.  The notations ${\mathbb N}$, ${\mathbb Z}$, ${\mathbb Q}$, ${\mathbb R}$, and ${\mathbb C}$ denote respectively the natural numbers, the integers, the rational numbers, the real numbers, and the complex numbers.  The complex projective space of complex dimension $n$ is denoted by ${\mathbb P}_n$.  The open disk in ${\mathbb C}$ of radius $r$ centered at the origin is denoted by $\Delta_r$.  When $r=1$, $\Delta_r$ is also denoted by $\Delta$.  We also use the notation $\Delta$ for the Laplacian where its meaning will be explicitly spelled out or clear from the context.

\bigbreak\noindent(0.5) {\it Remark on the Application of Arithmetic Technique Involving Heights to Analytic Objects Not Defined Over the Field of Algebraic Numbers.}  An arithmetic technique involving heights which was introduced by Gelfond and Schneider for the solution of Hilbert's seventh problem is used in Part I.  The key point in the use of heights is a conclusion about a positive lower bound to guarantee the nonvanishing of a limit. This lower bound is obtained from the self-evident statement that the norm of an integer over a number field which is nonzero is bounded from below by $1$.   In analysis it is difficult to obtain a positive lower bound to guarantee the nonvanishing of a limit, for which, for example, one of the available techniques is Harnack's inequality.  Here the arithmetic technique involving heights is assuming this r\^ole of producing a positive lower bound to guarantee the nonvanishing of a limit.  Since the application of this arithmetic technique involving heights requires the analytic objects under consideration to be defined over the field $\overline{\mathbb Q}$ of all algebraic numbers, we would like to remark how the general case is handled which does not satisfy this condition of definition over $\overline{\mathbb Q}$.

\medbreak For the situations which we consider we have a complex algebraic manifold $X$ defined over ${\mathbb C}$ by a finite number of homogeneous polynomials and we also have holomorphic vector bundles or other analytic objects ${\mathfrak F}$ which are defined algebraically over ${\mathbb C}$.  We can choose a finite number of complex numbers $c_1,\cdots,c_k$ which are algebraically independent over $\mathbb Q$ such that all these analytic objects ${\mathfrak F}$ under consideration are defined algebraically over a field $K$ which is a finite extension of the purely transcendental extension field ${\mathbb Q}\left(c_1,\cdots,c_k\right)$ of the rational number field ${\mathbb Q}$.   In this note our problem where the arithmetic technique involving heights is used deals with the existence or nonexistence of certain analytic objects.  We choose generic elements $\gamma_1,\cdots,\gamma_k$ in $\overline{\mathbb Q}$ and introduce another set of analytic objects $\tilde X$ and $\tilde{\mathfrak F}$ which are obtained by replacing $c_\ell$ by $\gamma_\ell$ in the coefficients of equations defining them.  When our conclusions hold for the set of analytic objects $\tilde X$ and $\tilde{\mathfrak F}$ for a generic choice of elements $\gamma_1,\cdots,\gamma_k$ in $\overline{\mathbb Q}$, the same conclusions hold for our original set of analytic objects $X$ and ${\mathfrak F}$.  In the sketch of our proof, when we have to use the arithmetic technique involving heights this way of handling the general case of analytic objects which may not be defined over $\overline{\mathbb Q}$ is implicitly assumed without any further explicit mention.

\bigbreak\noindent{\bf Table of Contents}

\bigbreak\noindent PART I. The Case of Zero Numerical Kodaira Dimension
\begin{itemize}
\item[\S1.] Curvature Current and Dichotomy
\item[\S2.] Gelfond-Schneider's Technique of Algebraic Values of Solutions of Algebraically Defined Differential Equations
\item[\S3.] Final Step of the Case of Zero Numerical Kodaira Dimension
\end{itemize}
\noindent PART II. The Case of General Numerical Kodaira Dimension Under the Assumption of Coincidence of Numerical Trivial Foliation and Fibration
\begin{itemize}
\item[\S4.] Numerically Trivial Foliations and Fibrations for Canonical Line Bundle
\item[\S5.] Curvature of Zeroth Direct Image of Relative Canonical and Pluricanonical Bundle
\item[\S6.] Strict Positivity of Direct Image of Relative Pluricanonical Bundle Along Numerically Trivial Fibers in the Base of Numerically Trivial Fibration
\end{itemize}
\noindent PART III. Coincidence of Numerically Trivial Foliation and Fibration for Canonical Bundle
\begin{itemize}
\item[\S7.] Technique of Nevanlinna's First Main Theorem for Proof of Compactness of Leaves of Foliation
\end{itemize}

\bigbreak

\bigbreak\noindent{\bf PART I. The Case of Zero Numerical Kodaira Dimension.}

\bigbreak\noindent{\bf\S1. Curvature Current and Dichotomy}

\bigbreak A general method of using the dichotomy in the decomposition of the curvature current to produce holomorphic sections of holomorphic line bundles was introduced in [Siu2006, \S6] for the analytic proof of the finite generation of the canonical ring [Siu2006, Siu2008, Siu2009].  The advantage of using curvature currents is that in the context of curvature currents one can use the roots of line bundles and the limits of flat twistings of multiples of line bundles with the desired amount of effective control.

\bigbreak\noindent(1.1) {\it Curvature Current of the Canonical Line Bundle.}  Let $L$ be a holomorphic line bundle over a compact complex algebraic manifold $X$ of complex dimension $n$ such that $\Gamma\left(X,m_0L\right)$ is nonzero for some positive integer $m_0$.  For any positive integer $m$ let
$$
s_1^{(m,L)},\cdots,s_{q_{m,L}}^{(m,L)}\in\Gamma\left(X,mL\right)
$$
be a ${\mathbb C}$-basis of $\Gamma\left(X,mL\right)$.  Let
$$
\Phi_L=\sum_{m=1}^\infty\varepsilon_{m,L}\sum_{j=1}^{q_{m,L}}\left|s^{(m,L)}_j\right|^{\frac{2}{m}},
$$
where $\varepsilon_{m,L}$ is a sequence of positive numbers going
decreasingly to $0$ fast enough to guarantee the convergence of the
infinite series.  Let
$$
\Theta_L=\frac{\sqrt{-1}}{2\pi}\,\partial\bar\partial\log\Phi_L.
$$
The reciprocal $\frac{1}{\Phi_L}$ of $\Phi_L$ defines a (possibly singular) metric along the fibers of $L$ whose curvature current $\Theta_L$ is a closed positive $(1,1)$-current on $X$.  The metric $\frac{1}{\Phi_L}$ as defined depends on the choice of $s_1^{(m,L)},\cdots,s_{q_{m,L}}^{(m,L)}$.

\medbreak Suppose the numerical Kodaira dimension $\kappa_{\rm num}(X)$ of $X$ is $\geq 0$.  We have $\dim_{\mathbb C}\Gamma\left(X, m_\nu K_X+A\right)\geq 1$ for some sequence of positive integers $m_\nu\to\infty$ and we can define the metric $$\frac{1}{\Phi_{m_\nu K_X+A}}$$ of $m_\nu K_X+A$ and its curvature current $\Theta_{m_\nu K_X+A}$ for $\nu\in{\mathbb N}$.  We can now define a curvature current $\Xi_{K_X}$ of the canonical line bundle $K_X$ of $X$ as the weak limit, in the space of closed positive $(1,1)$-currents on $X$, of
$$\frac{1}{m_{\nu_k}}\,\Theta_{m_{\nu_k}K_X+A}$$
as $k\to\infty$ for some subsequence $m_{\nu_k}$ of the sequence $m_\nu$.  It is possible to get such a weak limit because the total mass of
$$\frac{1}{m_\nu}\,\Theta_{m_{\nu}K_X+A}$$
with respect to a strictly positive smooth curvature form $\omega_A$ of the ample line bundle $A$ is equal to
$$\frac{1}{m_{\nu}}\left(m_{\nu}K_X+A\right)\cdot A^{n-1}$$
which is bounded uniformly in $\nu$.  The curvature current $\Xi_{K_X}$ of the canonical line bundle $K_X$ of $X$ depends on many choices.

\bigbreak\noindent(1.2) {\it Second Case of Dichotomy}.  The curvature current $\Xi_{K_X}$ defined in (1.1) admits the following decomposition as a closed positive $(1,1)$-current
$$
\Xi_{K_X}=\sum_{j=1}^J\gamma_j\left[V_j\right]+R,
$$
where $V_j$ is an irreducible nonsingular hypersurface of $X$ and $\gamma_j>0$ and the Lelong number of $R$ is $0$ outside a countable union of subvarieties of codimension $\geq 2$ in $X$ [Siu1974].  We have the dichotomy into the following two cases.
\begin{itemize}
\item[(i)] Either $J=\infty$ or $R\not=0$.
\item[(ii)] $J<\infty$ and $R=0$.
\end{itemize}
We claim that under the assumption
$$\sup_{m\in{\mathbb N}}\dim_{\mathbb C}\Gamma\left(X, mK_X+kA\right)\leq C_k<\infty\quad{\rm for\ all\ }k\in{\mathbb N}\leqno{(1.2.1)}$$
we must have the second case of the dichotomy with $J<\infty$ and $R=0$.  Choose a positive integer $\ell$ such that $\ell A$ is very ample.  If we are in the first case of the dichotomy with either $J=\infty$ or $R\not=0$, by the argument given in [Siu2006, \S6] (and also [Siu2009, Proposition (5.4)]) of using the curve which is the intersection of $n-1$ divisors of generic sections of $\ell A$, we conclude that
$$
\sup_{m\in{\mathbb N}}\dim_{\mathbb C}\Gamma\left(X,mK_X+\ell\left(n+2\right)A\right)=\infty
$$
which contradicts the assumption (1.2.1) and thus finishes the verification of the claim.  We now discuss how the rationality of the coefficients $\gamma_j$ (in possibly some other representation of $\Xi_{K_X}$) follows from a rather trivial linear algebra argument of linear independence of cohomology classes.

\bigbreak\noindent(1.3) {\it Rationality of Coefficients in Some Representation of Curvature Current as Element of the (1,1)-Cohomology Group.}  We consider $\Xi_{K_X}$ as an element of $H^{1,1}(X)$, which means that we use the same notation to denote the element of $H^{1,1}(X)$ defined by the closed positive $(1,1)$-current $\Xi_{K_X}$.  Here $H^{1,1}(X)$ is the Hodge group which is isomorphic to the set of harmonic $(1,1)$-forms on $X$ with respect to some chosen K\"ahler metric of $X$.  We can regard the equation
$$
\Xi_{K_X}=\sum_{j=1}^J\gamma_j\left[V_j\right]\leqno{(1.3.1)}
$$
also as an equation for elements of $H^{1,1}(X)$ with $J<\infty$.  We claim that there is another representation
$$
\Xi_{K_X}=\sum_{j=1}^J\hat\gamma_j V_j\leqno{(1.3.2)}
$$
as elements of $H^{1,1}(X)$ such that each $\hat\gamma_j$ is a positive rational number for $1\leq j\leq J$.

\medbreak The claim is verified as follows by using a rather trivial linear algebra argument of linear independence of cohomology classes.  By relabeling $V_1,\cdots,V_J$ we can assume without loss of generality that $V_1,\cdots,V_{J_0}$ are ${\mathbb R}$-linearly independent as elements of $H^{1,1}(X)$ and we can write
$$
V_k=\sum_{j=1}^{J_0}\beta_{k,j}V_j\quad{\rm for\ \ }J_0+1\leq k\leq J
$$
as elements in $H^{1,1}(X)$ for some $\beta_{k,j}\in{\mathbb Q}$.  In the space of closed positive $(1,1)$-currents we have
$$
\Xi_{K_X}=\sum_{j=1}^J\gamma_j V_j=\sum_{j=1}^{J_0}\gamma^*_jV_j+
\sum_{k=J_0+1}^J\gamma_k\left(V_k-\sum_{j=1}^{J_0}\beta_{kj}V_j\right),
$$
where
$$
\gamma_j^*=\gamma_j+\sum_{k=J_0+1}^J\gamma_k\beta_{kj}.
$$
As elements of $H^{1,1}(X)$ we have
$$
\Xi_{K_X}=\sum_{j=1}^{J_0}\gamma_j^* V_j.
$$
Since the elements of $H^{1,1}(X)$ induced by $\Xi_{K_X}$ and $V_1,\cdots,V_{J_0}$ all come
from $H^2\left(X,{\mathbb Q}\right)$ and since the elements of $H^{1,1}(X)$ induced by $V_1,\cdots,V_{J_0}$ are linearly independent, it follows that the numbers $\gamma_k^*$ are all rational for $J_0+1\leq k\leq J$.  For any rational numbers $\gamma_k^\prime$ for $J_0+1\leq l\leq J$, in the space of closed positive $(1,1)$-currents we have
$$
\displaylines{\Xi_{K_X}=\sum_{j=1}^{J_0}\gamma^*_jV_j+
\sum_{k=J_0+1}^J\gamma_k\left(V_k-\sum_{j=1}^{J_0}\beta_{kj}V_j\right)\cr=\sum_{j=1}^{J_0}\gamma^*_jV_j+
\sum_{k=J_0+1}^J\left(\gamma_k-\gamma_k^\prime\right)\left(V_k-\sum_{j=1}^{J_0}\beta_{kj}V_j\right)+
\sum_{k=J_0+1}^J \gamma_k^\prime\left(V_k-\sum_{j=1}^{J_0}\beta_{kj}V_j\right)\cr=\sum_{j=1}^J\hat\gamma_j V_j
+\sum_{k=J_0+1}^J\left(\gamma_k-\gamma^\prime_k\right)\left(V_k-\sum_{j=1}^{J_0}\beta_{kj}V_j\right),}
$$
where $\hat\gamma_j=\gamma^*_j-
\sum_{k=J_0+1}^J\gamma^\prime_k\beta_{kj}$ for $1\leq j\leq J_0$ and $\hat\gamma_j=\gamma^\prime_j$ for $J_0+1\leq j\leq J$ are all rational.  For $\gamma_k^\prime$ sufficiently close to $\gamma_k$ for $J_0+1\leq k\leq J$, the numbers
$$\hat\gamma_j=\gamma^*_j-
\sum_{k=J_0+1}^J\gamma^\prime_k\beta_{kj}=\gamma_j+
\sum_{k=J_0+1}^J\left(\gamma_k-\gamma^\prime_k\right)\beta_{kj}
$$
for $1\leq j\leq J_0$ and $\hat\gamma_j=\gamma^\prime_j$ for $J_0+1\leq j\leq J$ are all positive.  Thus as elements of $H^{1,1}(X)$ we have
$$
\Xi_{K_X}=\sum_{j=1}^J\hat\gamma_j V_j
$$
with $\hat\gamma_j$ all positive and rational for $1\leq j\leq J$.  This finishes the verification of the claim.

\bigbreak\noindent(1.4) {\it Nontrivial Section of Flatly Twisted Pluricanonical Line Bundle.}  Choose a sufficiently divisible positive integer $\hat m$ such that $\hat m\hat\gamma_j$ is a positive integer for $1\leq j\leq J$.  Then $\hat K_X$ is equal to the line bundle $\sum_{j=1}^J\hat m\hat\gamma_j V_j$ plus a flat line bundle.  Let $s_{V_j}$ be the canonical section of the line bundle $V_j$ on $X$.  Then $$\prod_{j=1}^J\left(s_{V_j}\right)^{\hat m\hat\gamma_j}$$
is a non identically zero holomorphic section of $\hat m K_X+F$ over $X$ for some flat line bundle $F$ over $X$.

\bigbreak\noindent{\bf\S2. Gelfond-Schneider's Technique of Algebraic Values of Solutions of Algebraically Defined Differential Equations}

\bigbreak\noindent(2.1) {\it Bombieri's Higher Dimensional Formulation.}  Gelfond [Gelfond1934] and Schneider [Schneider1934]
independently introduced the technique of algebraic values of solutions of algebraically defined differential equations to solve Hilbert's seventh problem [Hilbert1900], which later was put into the setting of group-homomorphisms and generalized to the higher-dimensional case by Lang and Bombieri [Lang1962, Lang1965, Lang1966, Bombieri1970, Bombieri-Lang1970].  The technique is actually more transparent when formulated for the higher-dimensional case.  We give here Bombieri's higher-dimensional formulation and indicate why the arguments work, in order to prepare for its comparison with the use in Part III later of the related argument of using the First Main Theorem in Nevanlinna theory (see (7.6) below).

\bigbreak\noindent(2.1.1){\it Theorem [Bombieri1970]}. Let $f_1,\cdots,f_{d+1}$ be meromorphic functions on ${\mathbb C}^d$ of finite order $\leq\rho$ in the sense that each is the quotient of two entire functions of growth order $O\left(e^{|z|^{\rho+\varepsilon}}\right)$ for any $\varepsilon>0$ where $z=\left(z_1,\cdots,z_d\right)$ are the coordinates of ${\mathbb C}^d$ and $\left|z\right|^2=\sum_{j=1}^d\left|z_j\right|^2$.  Let $K$ be a finite extension field of ${\mathbb Q}$. Assume that $f_1,\cdots,f_{d+1}$ do not satisfy any polynomial equation with coefficient in $K$ and that the following differential equation

    $$\frac{\partial f_j}{\partial z_k}=R_{jk}\left(f_1,\cdots,f_{d+1}\right)\leqno{(2.1.1.1)}$$ holds for $1\leq j\leq d+1$ and $1\leq k\leq d$, where $R_{jk}$ is a rational function in $d+1$ variables with coefficients in $K$.
Let $S$ be the set of all $z\in{\mathbb C}^n$ such that $f_j\left(z\right)\in K$ for $1\leq j\leq d+1$.
Then $S$ is contained in the zero-set of a polynomial of degree $\leq d(d+1)\rho[K:{\mathbb Q}]+2d$ in $z_1,\cdots,z_d$.

\bigbreak\noindent(2.1.2) {\it Sketch of the Main Arguments for its Proof}.  We give a sketch of its proof here for the purpose of highlighting some key points which will be used for comparison with the related argument in Part III of using the First Main Theorem of Nevanlinna theory to prove the coincidence of the numerically trivial foliation and fibration for the canonical line bundle (see (7.6) below).  We will do this sketch only for the special case where $f_1,\cdots,f_{d+1}$ are all entire, because the general case of meromorphic functions on ${\mathbb C}^d$ can be reduced to the case of entire functions on ${\mathbb C}^d$ by simply multiplying the $d+1$ meromorphic functions by an appropriate entire function on ${\mathbb C}^d$.  Let $m\in{\mathbb N}$ and $S_m$ be a subset of $m$ distinct points of $S$.  Choose a  non-identically-zero $P\left(\xi_1,\cdots,\xi_{n+1}\right)$ polynomial of degree $J$ in the indeterminates $\xi_1,\cdots,\xi_{d+1}$ with coefficients in $K$ such that $P\left(f_1,\cdots,f_{n+1}\right)$ vanishes to order at least $L$ at every point of $S_m$.  The algebraic independence of $f_1,\cdots,f_{d+1}$ over $K$ is to guarantee that $P\left(f_1,\cdots,f_{d+1}\right)$ cannot be identically zero on ${\mathbb C}^d$.  The polynomial $P\left(\xi_1,\cdots,\xi_{n+1}\right)$ in the indeterminates $\xi_1,\cdots,\xi_{n+1}$ is chosen by solving a system of linear equations with its coefficients $a_{j_1,\cdots,j_{d+1}}$ in $K$ for $0\leq j_\nu\leq J$ as unknowns where the linear equations are given by the vanishing of $P\left(f_1,\cdots,f_{d+1}\right)$ to order at least $J$ at every point of $S_m$.

\medbreak We also want appropriate height estimates for the coefficients of the polynomial $P\left(\xi_1,\cdots,\xi_{n+1}\right)$.  The height of $\frac{p}{q}$ with $p,q\in{\mathbb Z}$ relatively prime is the maximum of $|p|$ and $|q|$ and
its definition can be naturally generalized to elements of $K$.  Instead of using the terminology of heights we follow Bombieri's terminology in [Bombieri1970, p.278].  For $\alpha\in K$ let the denominator ${\rm den}(\alpha)$ be the smallest positive natural integer $d$ such that $d\alpha$ is an integer in $K$.  Let $\left\|\alpha\right\|=\max_\sigma\left|\sigma\alpha\right|$ and ${\rm size}(\alpha)=\max_\sigma\left(\log{\rm den}(\alpha), \log\left|\sigma\alpha\right|\right)$, where $\sigma$ ranges over all embeddings of $K$ into ${\mathbb C}$.  The solution of the linear equations to get $a_{j_1,\cdots,j_{d+1}}$ in $K$ with estimates of $\left\|a_{j_1,\cdots,j_{d+1}}\right\|$ and ${\rm den}\left(a_{j_1,\cdots,j_{d+1}}\right)$ is obtained by Siegel's lemma of the pigeon-hole principle [Bombieri1970, p.279, Lemma 2]. The differential equation $(2.1.1.1)$ is to guarantee that the coefficients of the linear equations are in $K$ and also to guarantee that we have the height estimate ${\rm size}\left(a_{j_1,\cdots,j_{d+1}}\right)=O\left(L\right)$ when $J^{d+1}$ is the integral part $\left\lfloor mL^d\log L\right\rfloor$ of $mL^d\log L$.

\medbreak For a closed positive $(1,1)$-current $\Xi$ on an open subset $\Omega$ of ${\mathbb C}^d$ let $\left\|\Xi\right\|$ be the mass measure of $\Xi$ which is the $(d,d)$-current defined by
$$\left\|\Xi\right\|=\frac{1}{(d-1)!}\,\Xi\wedge\left(\frac{\sqrt{-1}}{2}\sum_{j=1}^d dz_j\wedge d\overline{z_j}\right)^{d-1}.$$
For a subset $E$ of $\Omega$ we denote by $\left\|\Xi\right\|(E)$ the value of the measure $\left\|\Xi\right\|$ at the set $E$.
Let $s$ be the minimum vanishing order of $P\left(f_1,\cdots,f_{n+1}\right)$ at points of $S_m$ which is achieved at some point $Q$ of $S_m$. The key point of Bombieri's higher-dimensional version of the technique of Gelfond-Schneider is that, as $J\to\infty$ and $S_m$ exhausts $S$, the sequence of closed positive $(1,1)$-current
$$
\Theta=\frac{1}{s}\frac{\sqrt{-1}}{2\pi}\,\partial\bar\partial\log\left|P\left(f_1,\cdots,f_{n+1}\right)\right|^2
$$
satisfies the following uniform bound for its mass.

\medbreak\noindent(2.1.2.1) The mass measure $\left\|\Theta\right\|$ of $\Theta$ satisfies
$$
\left\|\Theta\right\|\left(B_r\right)\leq\frac{1}{{\rm Vol\,}_{d-1}}\left(\left(d+1\right)\rho+o(1)\right)\left[K:{\mathbb Q}\right]\quad{\rm for\ fixed\ }\ r.
$$
for fixed $r>0$, where $B_r$ is the ball $B_r$ of ${\mathbb C}^d$ of radius $r$ centered at the origin and ${\rm Vol\,}_{d-1}$ is the
the volume of the unit ball in ${\mathbb C}^{d-1}$.

\medbreak This estimate is obtained from the higher-dimensional version of Nevanlinna's technique of integration by parts twice.  In Nevanlinna's original paper [Nevanlinna1925], in order to derive his First Main Theorem [Nevanlinna1925, p.18, Erster Hauptsatz] he used the formula [Nevanlinna1925, p.5, $({\rm A})^{\prime\prime}$]
$$
\log|f(x)|=\frac{1}{2\pi}\int_{\partial G}\log|f(\xi)|\frac{\partial g(\xi,x)}{\partial n}\,ds+\sum_G g\left(x,b_\nu\right)-\sum_G g\left(x,a_\mu\right),
$$
where $f(x)$ is meromorphic on $G=\left\{\,x\in{\mathbb C}\,\big|\,r_1<|x|<r_2\,\right\}$ with zeroes $\left\{a_\mu\right\}$ and poles $\left\{b_\nu\right\}$ and $g$ is the Green's function of $G$.  This formula specializes to the special case of the Poisson-Jensen formula [Nevanlinna1925, p.5, (B)] that
$$
\log\left|f(0)\right|=\frac{1}{2\pi}\int_0^{2\pi}\log\left|f\left(re^{i\theta}\right)\right|d\theta-\sum_{\left|a_\mu\right|<r}\log\frac{r}{\left|a_\mu\right|}
+\sum_{\left|b_\nu\right|<r}\log\frac{r}{\left|b_\mu\right|},\leqno{(2.1.2.2)}
$$
where $f$ is a meromorphic function on ${\mathbb C}$ with zero-set $\left\{a_\mu\right\}$ and pole-set $\left\{b_\nu\right\}$ and has no zeroes and no poles at $0$ and on $\left\{\,x\in{\mathbb C}\,\big|\,\left|x\right|=r\,\right\}$.  Note that, when $f$ is entire, though the mass of the $(1,1)$-current $\frac{\sqrt{-1}}{2\pi}\partial\bar\partial\log\left|f\right|$ on the disk $\Delta_r$ of radius $r$ in ${\mathbb C}$ centered at $0$ is the number $n(r)$ of points $a_\mu$ in $\Delta_r$, the expression $\sum_{\left|a_\mu\right|<r}\log\frac{r}{\left|a_\mu\right|}$ is used in the above formula instead, because Green's function $g$ of $G$ involves the logarithmic function.  What is used in the technique of Gelfond-Schneider is that
$$
\sum_{\left|a_\mu\right|<r}\log\frac{r}{\left|a_\mu\right|}\geq \left(\log\frac{r}{r_0}\right)n\left(r_0\right)
$$
for $0<r_0<r$ so that it is possible to extract from it the factor $\log\frac{r}{r_0}$ which is of the same order as $\log r$ as $r\to\infty$ when we focus on the the lower bound of the number $n\left(r_0\right)$ of points $a_\mu$ in $\Delta_{r_0}$ for fixed $r_0$.
The following is the higher-dimensional version of (2.1.2.1) [Bombieri1970, p.273, Proposition 4] whose term for the mass measure of the closed positive $(1,1)$-current for the zero divisor contains also a logarithmic factor.

\medbreak\noindent(2.1.2.3) If $F$ is a holomorphic function on some open neighborhood of the topological closure of the open ball $B_R$ in ${\mathbb C}^d$ of radius $R$ centered at $0$, then for $0<\tau<1$ and $w_0\in B_{\frac{\tau R}{6d}}$
$$
\log\left|F(w_0)\right|\leq\max_{w\in\,\partial B_R}\log\left|F(w)\right|-(1-\tau)\left\|T\right\|\left(B_{\left|w_0\right|}\right)\log\frac{\frac{\,\tau R\,}{6d}}{\,\left|w_0\right|\,},
$$
where $T$ is the closed positive $(1,1)$-current $\frac{\sqrt{-1}}{2\pi}\partial\bar\partial\log\left|F\right|^2$ for the zero divisor of $F$ whose mass measure is denoted by $\left\|T\right\|$.

\medbreak From (2.1.2.3) an upper bound for $\left\|T\right\|\left(B_{\left|w_0\right|}\right)$ needs to come from a lower bound of $\log\left|F(w_0)\right|$ and an upper bound of $\max_{w\in\,\partial B_R}\log\left|F(w)\right|$.  By applying (2.1.2.3) to $F=P\left(f_1,\cdots,f_{d+1}\right)$ with $w_0$ averaged over $\partial B_r$ and with $R=s^\kappa$ for any $0<\kappa<\frac{1}{(d+1)\rho}$, we get
$$
(1-\tau)\kappa\,s\log s\,\left\|\Theta\right\|\left(B_r\right)\leq\left[K:{\mathbb Q}\right]\left(1+o(1)\right)s\log s\leqno{(2.1.2.4)}
$$
for any fixed $r>0$ with $Q\in B_{\frac{r}{2}}$ and for any $0<\tau<1$ as $L\to\infty$ and $S_m$ exhausting $S$ for the following reason.

\medbreak The factor $\kappa\log s$ on the left-hand side of (2.1.2.4) comes from $\log\frac{\frac{\,\tau R\,}{6d}}{\,\left|w_0\right|\,}$ with $R=s^\kappa$.  The factor $s$ on the left-hand side of (2.1.2.4) comes from the fact that $\Theta$ is obtained by dividing by $s$ the closed positive $(1,1)$-current for the zero divisor of $F$.  In the application of (2.1.2.3) the contribution from $\max_{w\in\,\partial B_R}\log\left|F(w)\right|$ is no more than $o\left(s\right)$, because the upper bound of $\log\left|F\right|$ on $B_R$ is no more than $$O\left(JR^{\rho+\varepsilon}\right)=O\left(\left(s^d\log s\right)^{\frac{1}{d+1}}s^{\kappa\left(\rho+\varepsilon\right)}\right)=o(s)$$ when $\varepsilon>0$ is chosen with $\kappa<\frac{1}{(d+1)\left(\rho+\varepsilon\right)}$.

\medbreak By applying the Cauchy integral formula for some nonzero partial derivative $\left(D^sF\right)(Q)$ of $F$ of order $s$ at $Q$ in terms of the value of $F$ on $\partial B_r$, one concludes that the average of $\log\left|F(w_0)\right|$ over $w\in\partial B_r$ is bounded from below by the sum of $-\left(1+o(1)\right)s\log s$ and $\log\left|\left(D^sF\right)(Q)\right|$, because the Cauchy integral formula for the $s$-th derivatives contains the factor $\frac{s!}{2\pi i}$ whose logarithm can be estimated by $s\log s\left(1+o(1)\right)$.  The lower bound $-\left(\left[K:{\mathbb Q}\right]-1\right)\left(1+o(1)\right)s\log s$ for $\log\left|\left(D^sF\right)(Q)\right|$ comes from differentiating (2.1.1.1) $s-1$ times and from the use of the chain rule, because we have
${\rm size}\left(a_{j_1,\cdots,j_{d+1}}\right)=O\left(L\right)$
and because there is the factor $\left[K:{\mathbb Q}\right]-1$ from the relation among ${\rm size}(\alpha)$, ${\rm den}(\alpha)$ and $\left|\sigma\alpha\right|$ for $\alpha\in K$ and an embedding $\sigma$ of $K$ into ${\mathbb C}$ [Bombieri1970, p.278].

\medbreak Since the Lelong number of the weak limit of $\Theta$ is at least $1$ at every point of $S$, the conclusion of Bombieri's higher-dimensional formulation now follows from the usual $L^2$ estimates of $\bar\partial$.

\bigbreak\noindent(2.1.2.5) {\it Gelfond-Schneider's Solution of Hilbert's 7th Problem}. The solution of Hilbert's seventh problem by Gelfond and Schneider [Gelfond1934, Schneider1934] is a special case of (2.1.1) with $f_1(z)=e^z$ and $f_2(z)=e^{\beta z}$ and $K={\mathbb Q}\left[\alpha, \alpha^\beta\right]$ so that if $\alpha^\beta$ is not transcendental for some algebraic $\alpha, \beta$ with $\alpha\not=0,1$ and $\beta$ irrational, then $S$ contains ${\mathbb Z}\log\alpha$, contradicting (2.1.1).

\bigbreak\noindent(2.1.3) {\it Formulation for Group-Homomorphism and Torsion Points}.  In the technique of Gelfond-Schneider the system of differential equations can be replaced by the differential of the group law when we have a holomorphic group-homomorphism from ${\mathbb C}^n$ to an algebraic group instead of $n+1$ meromorphic functions of finite order on ${\mathbb C}^n$ which are algebraically independent over the number field $K$.  This kind of formulation was introduced by Lang [Lang1962, Lang1965, Lang1966].  We need the following statement from this kind of formulation.

\medbreak\noindent(2.1.3.1)  Let $\varphi:{\mathbb C}^n\to G$ be a holomorphic map into a group $G$ which has an algebraic group structure over a number field $K$.  Assume that $\varphi$ is a group homomorphism.  Let $\vec{e_1},\cdots,\vec{e_n}$ be a ${\mathbb C}$-linear basis of ${\mathbb C}^n$ such that $\varphi\left(2\pi\sqrt{-1}\,\vec{e_j}\right)$ is a $K$-point of $G$ for $1\leq j\leq n$.  If the differential $d\varphi$ with respect to the ${\mathbb C}$-linear basis $\vec{e_1},\cdots,\vec{e_n}$ of ${\mathbb C}^n$ is algebraic and defined over $K$, then there is a finite subgroup $E$ of $G$ such that $\varphi\left(2\pi\sqrt{-1}\,\vec{e_j}\right)$ is in $E$ for $1\leq j\leq n$.

\medbreak\noindent(2.1.3.2) Here is its proof.  Without loss of generality we can assume that $\vec{e_1},\cdots,\vec{e_n}$ are the standard basis elements of ${\mathbb C}^n$.  The group $G$ can be embedded into some ${\mathbb P}_N$ with the embedding defined over $K$.  Let $f_1,\cdots,f_N$ be the $N$ meromorphic functions on ${\mathbb C}^n$ which are the pullbacks via $\varphi$ of the inhomogeneous coordinates of ${\mathbb P}_N$.  Let $f_{N+j}=e^{z_j}$ for $1\leq j\leq n$, where $z_1,\cdots,z_n$ are the $n$ coordinate functions of ${\mathbb C}^n$.  Let $\psi:{\mathbb C}^n\to  G\times{\mathbb C}^n$ be defined by $\varphi$ and the identity map of ${\mathbb C}^n$.   From the differential $d\psi$ and the group law of $G$ over $K$ and from the fact that the differential $d\varphi$ at the origin of ${\mathbb C}^n$ is defined over $K$ we get a system of differential equations on ${\mathbb C}^n$ whose unknowns are $N+n$ meromorphic functions $f_1,\cdots,f_{N+n}$ and which expresses any first-order partial derivative of each of $f_1,\cdots,f_{N+n}$ as a rational function of $f_1,\cdots,f_{N+n}$ with coefficients in $K$.  Note that since the group law of $G$ is defined over $K$, by the chain rule in order for the differential equations to be defined over $K$ it suffices to have $d\varphi$ defined over $K$ at one point of ${\mathbb C}^n$.

\medbreak The $N$ meromorphic functions $f_1,\cdots,f_N$ are of finite type, because $\varphi$ is a group homomorphism from which we can estimate the growth of meromorphic functions $f_1,\cdots,f_N$ on the ball of radius $r$ in ${\mathbb C}^n$ as $r\to\infty$ by using the group law of $G$.  Let $\Gamma$ be the discrete subgroup $\sum_{j=1}^n{\mathbb Z}\left(2\pi\sqrt{-1}\,\vec{e_j}\right)$ of ${\mathbb C}^d$.  Since $\Gamma$ is not contained in any algebraic hypersurface of ${\mathbb C}^n$ of finite degree and since $\psi\left(\Gamma\right)$ is contained in the algebraic points
of $G\times{\mathbb C}^n$, it follows from $(2.1.1)$ that the $N+n$ meromorphic functions $f_1,\cdots,f_{N+n}$ form a field of transcendence degree $\leq n$ over $K$.

\medbreak Let $\pi:{\mathbb C}^n\to{\mathbb C}^n\left/\Gamma\right.$ be the natural projection and let $\psi:{\mathbb C}^n\to G\times\left({\mathbb C}^n\left/\Gamma\right.\right)$ be defined by $\varphi$ and $\pi$.  Since ${\mathbb C}^n\left/\Gamma\right.$ is biholomorphic to $\left({\mathbb C}-\left\{0\right\}\right)^n$ under the map $\left(z_1,\cdots,z_n\right)\mapsto\left(e^{z_1},\cdots,e^{z_n}\right)$ and since the meromorphic functions
$f_1,\cdots,f_N,e^{z_1},\cdots,e^{z_n}$ form a field of transcendence degree $\leq n$ over $K$, it follows that the image $\psi\left({\mathbb C}^n\right)$ of ${\mathbb C}^n$ in $G\times\left({\mathbb C}^n\left/\Gamma\right.\right)$ is a subvariety of complex dimension $n$ in $G\times\left({\mathbb C}^n\left/\Gamma\right.\right)$.  Moreover, the restriction, to $\psi\left({\mathbb C}^n\right)$, of the projection $G\times\left({\mathbb C}^n\left/\Gamma\right.\right)\to{\mathbb C}^n\left/\Gamma\right.$ to its second factor makes $\psi\left({\mathbb C}^n\right)$ a finite-sheeted branch cover of ${\mathbb C}^n\left/\Gamma\right.$ whose fiber over the origin is projected by the natural projection $G\times{\mathbb C}^n\to G$ to a finite subset $E$ of $G$.  This means that $\varphi\left(\Gamma\right)$ is contained in the finite subset $E$ of $G$ and the subgroup $\varphi\left(\Gamma\right)$ of $G$ is a finite subgroup of $G$.  The proof is finished.

\bigbreak\noindent(2.1.3.3) {\it Remark on the Normalizing Factor.}  The normalizing factor $2\pi\sqrt{-1}$ used in (2.1.3.1) is the same one from the Cauchy integral formula or from the curvature form of a holomorphic line bundle with a Hermitian metric along its fibers.

\bigbreak\noindent(2.1.4) {\it Formulation for Group-Homomorphism and Coset of Subgroup}.  In our application of the technique of Gelfond-Schneider in the case of zero numerical Kodaira dimension the points and their images under a group-homomorphism are not known to be defined over $K$ as in the assumption of $(2.1.3)$, but only a coset and its image under a group-homomorphism is known to be defined over $K$ as subvarieties.  A modification is needed.  This modification was introduced by Simpson to enable him to draw conclusions about torsion translates of triple tori [Simpson1993, p.370, Corollary 3.5], which is used for the case of zero numerical Kodaira dimension.  We will not go into Simpson's setting to define what he means by torsion translates of triple tori.  Instead we will just extract from his work the simple argument which we need for our situation.  Here is the modification of (2.1.3.1) for the coset formulation.

\bigbreak\noindent(2.1.4.1)  Let $\varphi:{\mathbb C}^n\to G$ be a holomorphic map into a group $G$ which has an algebraic group structure over a number field $K$ such that $\varphi$ is a group homomorphism.  Let $N$ be a ${\mathbb C}$-linear subspace of ${\mathbb C}^n$ defined by a ${\mathbb C}$-basis over $K$ and let $v$ be a point of ${\mathbb C}^n$ such that the coset $N+v$ is defined by the vanishing of polynomials on ${\mathbb C}^n$ with coefficients in $K$.  Assume that $\varphi\left(N+v\right)$ is a subvariety of $G$ defined over $K$.  If the normalized differential $\frac{\sqrt{-1}}{2\pi}d\varphi$ at the origin of ${\mathbb C}^n$ is algebraic and defined over $K$, then there is an element $w$ of $N+v$ such that $\varphi(w)$ is a torsion point of $G$.

\bigbreak\noindent(2.1.4.2)  The proof is as follows.  In the application of (2.1.3.1) we replace $G$ by the quotient $G/\varphi(N)$ and replace ${\mathbb C}^n$ by the ${\mathbb C}$-line in ${\mathbb C}^n/N$ containing the image $\hat v$ of $v$ in ${\mathbb C}^n/N$ and replace $\vec{e_1}$ by $\hat v$.  Then from the application of (2.1.3.1) with $n=1$ we conclude that $\hat v$ is a torsion element of $G/\varphi(N)$ and as a consequence there is an element $w$ of $N+v$ such that $\varphi(w)$ is a torsion point of $G$.

\bigbreak\noindent(2.2) {\it Two Algebraic Descriptions of Flat Line Bundles.}  In order to apply (2.1.4.1) we need two different algebraic descriptions of flat line bundles on a compact complex manifold $X$ so that one which plays the r\^ole of the algebraic group $G$ is the integration of the other which plays the r\^ole of ${\mathbb C}^n$ and both are related by the holomorphic group-homomorphism $\varphi$.  For the algebraic description which plays the the r\^ole of ${\mathbb C}^n$ we use the set of all holomorphic $1$-forms with their algebraic structure for the description.  For the algebraic description which plays the the r\^ole of the algebraic $G$ we use transition functions with their algebraic structure for the description.

\bigbreak\noindent(2.2.1) Let $L$ be a flat line bundle over $X$ which is described by a finite open cover $\left\{U_j\right\}_{j\in J}$ of $X$ and transition functions $g_{jk}$ on $U_j\cap U_k$ such that each $g_{jk}$ is constant with $\left|g_{jk}\right|\equiv 1$ on $U_j\cap U_k$ and
there are local holomorphic frames $F_j$ of $L$ on $U_j$ ({\it i.e.,} each $F_j$ is a holomorphic section of $L$ on $U_j$ generating ${\mathcal O}(L)$ on $U_j$) with $F_j=g_{jk}F_k$ on $U_j\cap U_k$. For any section $s=\left\{s_j\right\}$ of $L$ on an open subset $G$ of $X$ with $s_j$ being a function on $G\cap U_j$ we can use $\left\|s\right\|^2=\left|s_j\right|^2$ which is independent of $j$ as a Hermitian metric along the fibers of $L$ whose curvature form is identically zero on $X$.

\medbreak For any holomorphic $1$-form $\omega$ on $X$ and for any complex number $\lambda\in{\mathbb C}$, we define a new complex structure $L_{\lambda,\omega}$ on the smooth line bundle $L$ as follows.  A local smooth section $s$ of $L_\lambda$ is defined to be holomorphic if and only if $\bar\partial s=\lambda s\,\bar\omega$.  In other words, the new differentiation in the $(0,1)$-direction for local smooth sections of $L_\lambda$ is $\bar\partial_\lambda=\bar\partial -\lambda\,\bar\omega$.

\bigbreak\noindent(2.2.2) We can determine the transition functions for $L_{\lambda,\omega}$ as follows.  By replacing the finite open cover $\left\{U_j\right\}_{j\in J}$ of $X$ by a refinement if necessary, we can assume without loss of generality that there exists a holomorphic function $f_j$ on $U_j$ such that $df_j=\omega$ on $U_j$ and $f_j-f_k=c_{kj}$ for some constant $c_{jk}$ on $U_j\cap U_k$.  Let $s_{j,\lambda}=e^{\lambda\overline{f_j}}F_j$ on $U_j$.  Then $\bar\partial s_{j,\lambda}=\lambda\,s_{j,\lambda}\bar\omega$ on $U_j$.  The collection $$\left\{s_{j,\lambda}\right\}_{j\in J}=\left\{e^{\lambda\overline{f_j}}F_j\right\}_{j\in J}$$ can be used as local holomorphic frames for the holomorphic line bundle $L_\lambda$ on $X$.  The transition functions $\left\{g_{jk,\lambda,\omega}\right\}_{j,k\in J}$ for $L_{\lambda,\omega}$ with respect to the cover $\left\{U_j\right\}_{j\in J}$ are given by
$$g_{jk,\lambda,\omega}=\frac{\ e^{\lambda\overline{f_j}}F_j\ }{\ e^{\lambda\overline{f_k}}F_k}=
e^{\lambda\left(\overline{f_j}-\overline{f_k}\right)}g_{jk}=e^{-\lambda\overline{c_{jk}}}g_{jk}
$$
on $U_j\cap U_k$ which is a constant.  One conclusion we can immediately get from this explicit expression for the transition function is that the line bundle $L_{\lambda,\omega}$ is holomorphic in the variable $\lambda\in{\mathbb C}$.

\bigbreak\noindent(2.2.3)  For reasons explained in Remark(0.4) we can assume without loss of generality that $X$ as well as all relevant analytic objects over $X$ which come into consideration in our discussion is assumed to be defined over $\overline{\mathbb Q}$.  We can describe the two algebraic structures for the set of flat line bundles over $X$.

\medbreak The first description is to give the flat line bundle $L$ with constant transition functions $g_{jk}$ (not necessarily of absolute-value $1$) the algebraic structure for the complex variables $g_{jk}$.  This is the same as giving the set of flat line bundles the complex structure of
${\rm Hom}\left(\pi_1(X),{\mathbb C}^*\right)={\rm Hom}\left(H_1\left(X,{\mathbb Z}\right),{\mathbb C}^*\right)$ defined by the generators of $\pi_1(X)$ (or
$H_1\left(X, {\mathbb Z}\right)$) and their relations.  We call this description the {\it algebraic structure from transition functions.}

\medbreak The second description is to give $L_{\lambda,\omega}$ the algebraic structure for the holomorphic line bundle $\omega$ and for the complex variable $\lambda$ when the line bundle $L$ which we start out with is the trivial line bundle.  We call this description the {\it algebraic structure from holomorphic $1$-forms}.

\medbreak As we have mentioned in (2.2) the algebraic structure from transition functions is obtained from the algebraic structure from holomorphic $1$-forms by using integration or the process of finding local primitives of $1$-forms.

\bigbreak\noindent(2.3) {\it Coset Defined by Lower Bound of Flatly Twisted Hodge Number.}  In order to apply the formulation of the technique of Gelfond-Schneider for group-homomorphisms and cosets of groups as given in (2.1.4.1), we would like to produce a coset in the second description of algebraic structures for the space of flat bundles.  For this purpose we have the following statement which is a special case of the general result of Simpson on torsion translates of triple tori [Simpson1993].

\bigbreak\noindent(2.3.1) Let $X$ be a compact complex algebraic manifold of complex dimension $n$ and $L$ is a flat line bundle over $X$.  Let $m$ be the complex dimension of the space $\Gamma\left(X,T_X^*\right)$ of all holomorphic $1$-forms on $X$.  Let $\omega$ be a holomorphic $1$-form on $X$.  Then for any nonzero complex number $\lambda\in{\mathbb C}^*$ the dimension of $\Gamma\left(X, L_{\lambda,\omega}+K_X\right)$ is independent of $\lambda\in{\mathbb C}^*$ (but may depend on $\omega$).  As a consequence, for any $k\in{\mathbb N}$ the set of all flat line bundles $F$ (as points of ${\mathbb C}^m$ and given the algebraic structure from holomorphic $1$-forms) which satisfy $\dim_{\mathbb C}\Gamma\left(X, F+K_X\right)\geq k$ is of the form $N+v$ for some ${\mathbb C}$-linear subspace $N$ of ${\mathbb C}^m$ and some $v\in{\mathbb C}^m$.

\bigbreak\noindent(2.3.2) For the following reason the above statement is a consequence of the Hodge decomposition with flat twisting.  Let $H^n_{\rm top}\left(X,F\right)$ be the topologically defined $n$-th cohomology group of $X$ with coefficients in the local system given by the flat line bundle $F$.  By standard arguments of topology and homological algebra the complex dimension of the topologically defined $H^n_{\rm top}\left(X,L_{\lambda,\omega}\right)$ is independent of $\lambda\in{\mathbb C}^*$ for any prescribed holomorphic $1$-form $\omega$ on $X$, because the transition function $g_{jk,\lambda,\omega}$ of $L_{\lambda,\omega}$ on $U_j\cap U_k$ is given by $e^{-\lambda\overline{c_{jk}}}g_{jk}$, where $g_{jk}$ is the transition function for $L$ on $U_j\cap U_k$ and $c_{kj}=f_j-f_k$ with $df_j=\omega$ on $U_j$.  Note that these arguments of topology and homological algebra work only for $\lambda\in{\mathbb C}^*$ and not for $\lambda\in{\mathbb C}$ when $\lambda\to 0$.

\medbreak By Hodge theory with flat twisting $L_{\lambda,\omega}$, we know that $\Gamma\left(X, L_{\lambda,\omega}+K_X\right)$ is a summand in a direct sum decomposition of $H^n_{\rm top}\left(X,L_{\lambda,\omega}\right)$ where the complex dimension of each summand is an upper semi-continuous function of $\lambda\in{\mathbb C}^*$.  Hence the dimension of $\Gamma\left(X, L_{\lambda,\omega}+K_X\right)$ is independent of $\lambda\in{\mathbb C}^*$ for any prescribed holomorphic $1$-form $\omega$ on $X$.  By linear algebra arguments, for any $k\in{\mathbb N}$ the set of all flat line bundles $F$ (as points of ${\mathbb C}^m$ and given the algebraic structure from holomorphic $1$-forms) which satisfy $\dim_{\mathbb C}\Gamma\left(X, F+K_X\right)\geq k$ is of the form $N+v$ for some ${\mathbb C}$-linear subspace $N$ of ${\mathbb C}^m$ and some $v\in{\mathbb C}^m$.

\bigbreak\noindent(2.3.3) {\it Remark on Dimension Jump of Cohomology for Flat Line Bundles Parametrized by Nonzero Complex Numbers at Zero Value of Parameter.}  Note that in (2.3.1) there is in general a jump of  $\dim_{\mathbb C}\Gamma\left(X, L_{\lambda,\omega}+K_X\right)$ when $\lambda\in{\mathbb C}^*$ approaches $0$, otherwise $\dim_{\mathbb C}\Gamma\left(X, F+K_X\right)$ would be independent of the flat line bundle $F$, which is not true.  This phenomenon is related to the fact that the natural metric of line bundle ${\mathcal L}_\omega$ over $X\times{\mathbb C}$ whose restriction to $X\times\left\{\lambda\right\}$ is $L_{\lambda,\omega}$ does not have a semipositive curvature form.  The reason is as follows. There is a natural Hermitian metric $e^{-\varphi_{\lambda,\omega}}$ along the fibers of $L_{\lambda,\omega}$ which on $U_j$ is given by $e^{-\varphi_{j,\lambda,\omega}}=e^{-2{\rm Re}\left(\lambda\overline{f_j}\right)}$ so that
$$
e^{-\varphi_{j,\lambda,\omega}}\left|g_{jk,\lambda,\omega}\right|^2=e^{-2{\rm Re}\left(\lambda\overline{f_j}\right)}\left|e^{\lambda\left(\overline{f_j}-\overline{f_k}\right)}g_{jk}\right|^2=
e^{-2{\rm Re}\left(\lambda\overline{f_k}\right)}=e^{-\varphi_{j,\lambda,\omega}}
$$
on $U_j\cap U_k$.

\medbreak A natural Hermitian metric $e^{-\varphi_{{\mathcal L}_\omega}}$ for the line bundle ${\mathcal L}_\omega$ on $X\times{\mathbb C}$ is given by
$$e^{-\varphi_{{\mathcal L}_\omega}}=e^{-\varphi_{j,\lambda,\omega}}=e^{-2{\rm Re}\left(\lambda\overline{f_j}\right)}$$
on $U_j\times{\mathbb C}$.  Its curvature current in general may not be semipositive, because that the function $2{\rm Re}\left(\lambda\overline{f_j}\right)$, though pluriharmonic in the variables along $X$ so that the curvature of $e^{-\varphi_\lambda}$ is zero on $X$ for a fixed $\lambda$, may {\it not} be pluriharmonic or plurisubharmonic {\it jointly} in both the variable $\lambda$ and the variables along $X$, because the function $\overline{f_j}$ is locally conjugate-holomorphic while the function $\lambda$ is holomorphic.

\medbreak If the curvature form of ${\mathcal L}_\omega$ were semipositive, by the extension theorem of Ohsawa-Takegoshi [Ohsawa-Takegoshi1987] (or the pluricanonical case in [Siu1998, Siu2002, Paun2007, Varolin2008]), $\dim_{\mathbb C}\Gamma\left(X, L_{\lambda,\omega}+K_X\right)$ (and even $\dim_{\mathbb C}\Gamma\left(X, L_{\lambda,\omega}+qK_X\right)$ for any fixed $q\in{\mathbb N}$) would have been independent of $\lambda\in{\mathbb C}$, which is not true when $\lambda\to 0$.  If we could replace $\lambda$ by $\overline{\lambda}$ or $\frac{1}{\ \overline{\lambda}\ }$ (which holds only when $|\lambda|=1$) or some other holomorphic function in $\overline{\lambda}$, the curvature form of ${\mathcal L}_\omega$ would be semipositive, which of course is impossible because we need to make the family of holomorphic line bundles $L_{\lambda,\omega}$ holomorphic in the variable $\lambda$.

\bigbreak\noindent(2.4) {\it Simpson's Result on Torsion Flat Bundle for Flatly Twisted Canonical Section.}  Let $X$ be a compact complex algebraic manifold of complex dimension $n$ and $G$ be the moduli space of all flat line bundles on $X$ given the algebraic structure from the transition functions.  For $\gamma\in G$ let $F_\gamma$ be the flat line bundle on $X$ associated to $\gamma$.  For $k\in{\mathbb N}$ let $V_k$ be the set of all $\gamma\in G$ such that $\dim_{\mathbb C}\Gamma\left(X, F_\gamma+K_X\right)\geq k$.  Then $V_k$ is an algebraic subvariety of $G$ from the general theory in algebraic geometry concerning the behavior of the dimension of the cohomology group of a varying fiber.

\medbreak By letting $L^{(0)}$ be the trivial line bundle on $X$ and setting $\lambda=1$, to a holomorphic $1$-form $\omega$ we can assign the flat line bundle $L^{(0)}_{1,\omega}$ to $\omega$.  Let $m$ be the complex dimension of the space of all holomorphic $1$-forms on $X$ so that the space of all holomorphic $1$-forms on $X$ can be identified with points of ${\mathbb C}^m$.  We now consider the map $\varphi:{\mathbb C}^m\to G$ which assigns to each holomorphic $1$-form $\omega$ on $X$ the element $\gamma$ of $G$ such that the flat line bundle $L^{(0)}_{1,\omega}$ is the same as $F_\gamma$.  Now the following special case of a result of Simpson [Simpson1993] is a consequence of (2.3.1) and (2.1.4.1).

\bigbreak\noindent(2.4.1) Let $X$ be a compact complex algebraic manifold and $k\in{\mathbb N}$.  If there exists a flat line bundle $F$ over $X$ such that $\dim_{\mathbb C}\Gamma\left(X, F+K_X\right)\geq k$, then there exists a flat line bundle $E$ on $X$ such that $pE$ is equal to the trivial line bundle on $X$ for some positive integer $p$ and $\dim_{\mathbb C}\Gamma\left(X, E+K_X\right)\geq k$.

\bigbreak\noindent(2.4.2) {\it Brieskorn's Use of Algebraic Values of Solutions of Algebraically Defined Differential Equations.}  We would like to remark that as mentioned by Simpson in [Simpson1993, p.362], Brieskron in [Brieskorn1070, p.114, (4), (5), (6) in Proof of Satz 4] already used Gelfond-Schneider's method of algebraic values of solutions of algebraically defined differential equations to draw conclusions on the monodromy of isolated singularities of hypersurfaces in arguments similar to Simpson's.  What is new in Simpson's technique is the use of the quotient of abelian variety with respect to an algebraically defined abelian subvariety in order to get an algebraic point in the quotient [Simpson1993 p.370, Proof of Theorem 3.3] as explained in (2.1.4.1) and (2.1.4.2).

\bigbreak\noindent{\bf\S3. Final Step of the Case of Zero Numerical Kodaira Dimension}

\bigbreak\noindent(3.1) {\it Use of Branched Cover to Construct Pluricanonical Section From Flatly Twisted Pluricanonical Section.}   For the case of a compact complex algebraic manifold $X$ with zero numerical Kodaira dimension, we have already got to the point of the existence of a nontrivial flatly twisted pluricanonical section on $X$.  What is left is to get rid of the flat twisting which we will do by lifting the flatly twisted pluricanonical section to a flatly twisted canonical section on the manifold $\hat X$ constructed from a branched cover over $X$ and then pushing forward to get a nontrivial pluricanonical section on $X$ after we invoke (2.4.1) to get a nontrivial pluricanonical section on $\hat X$.

\bigbreak\noindent(3.1.1) {\it Construction of Manifold from Branched Cover and Lifting of Flatly Twisted Pluricanonical Section to Flatly Twisted Canonical Section.}  Suppose $L$ is a flat line bundle over a compact complex algebraic manifold $X$ of zero numerical Kodaira dimension and suppose $m$ is a positive integer $\geq 2$.  Assume we have a nonzero element $s$ of $\Gamma\left(X, L+m K_X\right)$.  We would like to construct from it a nonzero element of $\Gamma\left(X, m^\prime K_X\right)$ for some positive integer $m^\prime$.  Let $N$ be the order of the torsion subgroup of $H^2\left(X,{\mathbb Z}\right)$ which is the same as the order of the torsion subgroup of $H_1\left(X,{\mathbb Z}\right)$ by the universal coefficient theorem.  The flat line bundle $NL$ belongs to the identity component of the real-analytic group of all flat line bundles over $X$.  As a result, we can write $NL$ as $mNL^\prime$ for some flat line bundle $L^\prime$ over $X$.  By replacing $s$ by $s^N$ and $L$ by $NL=mNL^\prime$ and then replacing $mN$ by $m$, we can assume without loss of generality that we have a nonzero element $s$ of $\Gamma\left(X, L+m K_X\right)$ such that $L=mL^\prime$ for some flat line bundle $L^\prime$ over $X$.

\medbreak Let $\tilde X$ be the graph of the multi-valued section $s^{\frac{1}{m}}$ of $L^\prime+K_X$ over $X$ and let $\pi:\tilde X\to X$ be the natural projection.  Then from $s^{\frac{1}{m}}$ we get a holomorphic section $s^\prime$ of $\pi^*\left(L^\prime+K_X\right)$ over $\tilde X$.  At regular point of the divisor of $s$ defined by $\zeta=0$ locally, an $m$-th root $w$ of $\zeta$ is given by $\zeta=w^m$ which yields $d\zeta=mw^{m-1}dw$ so that the Jacobian determinant of $\pi:\tilde X\to X$ is equal to $m-1$ times the divisor of $s^\prime$ at a regular point of the divisor of $s$.  Let $\hat X$ be the desingularization of (the normalization of) $\tilde X$ with projection $\hat\pi:\hat X\to X$ which is the composite of the desingularization projection $\hat X\to\tilde X$ and $\pi:\tilde X\to X$.  Then the canonical line bundle $K_{\hat X}$ is equal to $\hat\pi^*\left(K_X+(m-1)\left(K_X+L^\prime\right)\right)+D$ for some nonnegative divisor $D$ on $\hat X$ whose $\hat\pi$ image is contained in the zero-set of $s$.  Let $s_D$ be the canonical section of the line bundle over $\hat X$ associated to $D$.  Since $\hat\pi^*s$ is an element of $\Gamma\left(\hat X, \hat\pi^*\left(mK_X+mL^\prime\right)\right)$, it follows that $s_D\,\hat\pi^*s$ is a nonzero element of $$\Gamma\left(\hat X, \hat\pi^*\left(K_X+(m-1)\left(K_X+L^\prime\right)+L^\prime\right)+D\right)=\Gamma\left(\hat X, K_{\hat X}+\hat\pi^*L^\prime\right).$$

\bigbreak\noindent(3.1.2) {\it Removal of Flat Twisting and Pushing Forward to Get Pluricanonical Section.}  The complex dimension of $\Gamma\left(\hat X, K_{\hat X}+\hat\pi^*L^\prime\right)$ is $\geq 1$, because it contains the nonzero element $s_D\,\hat\pi^*s$.  By (2.4.1) there is a flat line bundle $L^{\#}$ on $\hat X$ such that the complex dimension of $\Gamma\left(\hat X, K_{\hat X}+L^{\#}\right)$ is $\geq 1$ and $\hat mL^{\#}$ is the trivial line bundle on $\hat X$ for some positive integer $\hat m$.  Hence there exists some
nonzero element $\hat s$ in $\Gamma\left(\hat X,\hat mK_{\hat X}\right)$ for some positive integer $\hat m$.

\medbreak By taking the product of the values of $\hat s$ on the $m$ sheets of $\hat X$ over $X$ outside the divisor of $s$, we get a non identically zero {\it meromorphic} section $\sigma$ of $m\hat mK_X$ over $X$ whose pole-set is contained in the divisor of $s$.  We claim that $\sigma$ must be holomorphic.  Assume the contrary and we are going to derive a contradiction.  Since the pole-set of $\sigma$ is contained in the divisor of $s$, it follows that for some sufficiently large $N$ the product $\sigma s^N$ is a {\it holomorphic} section of the line bundle $m\left(\hat m +N\right)K_X+NL$ over $X$ whose divisor is different from that of $s^{\hat m+N}$.   Since $L$ is flat, we can choose its transition functions to be constant of absolute value $1$ with respect to some local holomorphic frames for some finite cover of $X$ so that, for a local section $s$ of $L$ regarded as a local function $f$ by using such a local frame, we can use $|f|$ as the pointwise norm of the local section $s$.  With such a pointwise norm for $L$, we can define the metric
$$\frac{1}{\left(\left|\sigma s^N\right|^2+\left|s^{\hat m +N}\right|^2\right)^{\frac{1}{m\left(\hat m +N\right)}}}$$ for $K_X$ and use its curvature current $\hat\Theta$ to compute the growth rate of the dimension of $\Gamma\left(X, pK_X+(n+2)\hat A\right)$ as a function of $p$ for any very ample line bundle $\hat A$ over $X$.  Since the divisor of the product $\sigma s^N$ is different from that of $s^{\hat m+N}$, it follows that the curvature current $\hat\Theta$ is not in the second case of the dichotomy and the dimension of $\Gamma\left(X, pK_X+(n+2)\hat A\right)$ is no longer a bounded function of $p$.  This contradicts the assumption of zero numerical Kodaira dimension.  Thus we know that $\sigma$ is a non identically zero holomorphic section of $\hat m m K_X$ over $X$.  This concludes the confirmation of the abundance conjecture for the case of zero numerical Kodaira dimension.

\bigbreak\noindent(3.2) {\it Remarks.} (a) The assumption of zero numerical Kodaira dimension is essential to the simple argument of using a branched cover to construct pluricanonical sections from flatly twisted pluricanonical sections.

\medbreak\noindent(b) The posted paper of Campana Peternell, and Toma [Campana-Peternell-Toma2007] already used the technique of applying Simpson's result of (2.4.1) to a branched cover.  The only difference is that the argument in [Campana-Peternell-Toma2007] is algebraic geometric involving divisors and here the argument is analytical involving metrics and curvature.

\medbreak\noindent(c) Nero Budur informed me that his generalization of Simpson's result on torsion translates of triple tori [Simpson1993] in [Budur2009, Theorem 1.3] involving multiplier ideals can be used to offer an alternative to the argument of branched cover of (3.1.1).

\bigbreak\noindent{\bf PART II. The Case of General Numerical Kodaira Dimension Under the Assumption of Coincidence of Numerical Trivial Foliation and Fibration.}

\bigbreak\noindent{\bf\S4. Numerically Trivial Foliations and Fibrations for Canonical Line Bundle.}

\bigbreak For any holomorphic line bundle $L$ over a compact complex algebraic manifold $X$, since the introduction of Kodaira dimension and numerical Kodaira dimension and the later work of Iitaka [Iitaka1971], there have been quite a number of investigations about factoring out the directions along which the line bundle $L$ is trivial in various senses (for example, [Tsuji2000, Bauer-et-al2002, Takayama2003, Eckl2004a, Siu2008, Siu2009] to name some).  In this Part II we will introduce in a very brief manner the numerical trivial foliation [Eckl2004b] and the numerically trivial fibration [Tsuji2000, Bauer-et-al2002, Takayama2003, Eckl2004a] for a general holomorphic line bundle $L$. Then we specialize to the case of the canonical line bundle and confirm the abundance conjecture for general numerical Kodaira dimension under the additional assumption that it is known that the numerically trivial foliation for the canonical line bundle always coincides with the numerically trivial fibration for the canonical line bundle, leaving to Part III the verification that this additional assumption holds always automatically as a result of the application of a modification of the technique of Gelfond-Schneider [Gelfond1934, Schneider1934].

\bigbreak\noindent(4.1) {\it Definitions of Numerically Trivial Foliation and fibration.}  Recall that in (1.1) the metric $\frac{1}{\Phi_L}$ is defined for a holomorphic line bundle $L$ on a compact complex algebraic manifold $X$ of complex dimension $n$.  For an ample line bundle $A$ and the line bundle $L$, as was done for $L=K_X$ in (1.1), we can introduce ample-bundle-defined curvature current $\Xi_L$ as the weak limit, in space of closed positive $(1,1)$-currents on $X$, of
$$\frac{1}{m_\nu}\,\frac{\sqrt{-1}}{2\pi}\,\partial\bar\partial\log\Phi_{m_\nu L+A}$$
as $\nu\to\infty$ for some sequence $m_\nu\to\infty$ with $m_\nu\in{\mathbb N}$.
Let
$\Xi_L=\sum_{j\in J}c_j\left[V_j\right]+R$ be its canonical decomposition [Siu1974] where each $c_j>0$, each $V_j$ is an irreducible complex hypersurface of $X$, and $R$ is a closed positive $(1,1)$-current on $X$ whose Lelong number is $0$ outside at most a countable union $\cup_{\lambda\in\Lambda}W_\lambda$ of irreducible complex subvarieties $W_\lambda$ of complex codimension $\geq 2$ in $X$.

\medbreak For any irreducible complex curve $C$ we say that $L$ is numerically trivial on $C$ if $C$ is disjoint from $\cup_{\lambda\in\Lambda}W_\lambda$ and for an open subset $U$ of $X$ which intersects $C$ on which $R=\frac{\sqrt{-1}}{2\pi}\,\partial\bar\partial\psi$ for some plurisubharmonic function $\psi$, the restriction of $\psi$ to $C\cap U$ is harmonic (and hence pluriharmonic).  Likewise, for an irreducible locally closed complex curve $E$ of $X$ we say that $L$ is numerically trivial on $E$ if $E$ is disjoint from $\cup_{\lambda\in\Lambda}W_\lambda$ and for an open subset $U$ of $X$ which intersects $E$ on which $R=\frac{\sqrt{-1}}{2\pi}\,\partial\bar\partial\psi$ for some plurisubharmonic function $\psi$, the restriction of $\psi$ to $E\cap U$ is harmonic (and hence pluriharmonic).  Here by a locally closed complex curve $E$ in $X$ we mean that $E$ is a complex curve in some open subset of $X$.  Note that in this definition of numerical triviality of $L$ on an irreducible complex curve or locally closed complex curve the part $\sum_{j\in J}c_j\left[V_j\right]$ of $\Xi_L$ is not involved.  This non-involvement of $\sum_{j\in J}c_j\left[V_j\right]$ makes the definition of numerical triviality for an irreducible complex curve $C$ different from the condition that $L\cdot C=0$.

\medbreak The numerically trivial fibration is obtained by piecing together chain-wise intersecting numerically trivial irreducible complex curves to form the fibers of the fibration.  The numerically trivial fibration may be defined only outside a pluripolar subset $Z$ of $X$ so that it is a holomorphic map from $X-Z$ to some complex space $S$.  We may have to blow up $X$ in order to make the numerically trivial fibration given by a holomorphic from $X$ to some complex space $S$.  Such blow-ups do not matter much in our investigation, because pluricanonical sections are invariant under blowups.

\medbreak The numerically trivial foliation is obtained by putting together the tangent bundle of all the numerically trivial irreducible local complex curves to form an analytic subvariety of the tangent bundle of $X$.

\medbreak For more detailed descriptions see [Tsuji2000, Bauer-et-al2002, Takayama2003, Eckl2004a] for numerically trivial fibrations and see [Eckl2004b] for numerically trivial foliations.

\medbreak In this Part II we assume the coincidence of numerical trivial foliation and numerical trivial fibration for the canonical line bundle whose verification will be presented in Part III.  Here in this Part II present the proof of the abundance conjecture for the case of positive numerical Kodaira dimension under this additional assumption.  After applying blow-ups to $X$ if necessary, we now assume that there is a holomorphic map $\pi:X\to S$ from $X$ to a compact complex manifold $S$ such that the generic fibers of $X$ have zero numerical Kodaira dimension and the complex dimension of $S$ is equal to the numerical Kodaira dimension of $X$.

\medbreak Let $K_{X|S}$ be the relative canonical line bundle of $X$ whose restriction to a generic fiber is the canonical line bundle of the fiber.  The key point here is that from $K_X=K_{X|S}\left(\pi^*K_S\right)$ the zeroth direct image $R^0\pi_*\left(mK_{X|S}\right)$ of the relative pluricanonical line bundle $mK_{X|S}$ plays a r\^ole in the determination of the numerical Kodaira dimension of $X$.  Since the curvature of $R^0\pi_*\left(mK_{X|S}\right)$ is semipositive on $S$ (see (5.9)(a)), the problem is reduced to verifying that the curvature of $R^0\pi_*\left(mK_{X|S}\right)$ is strictly positive on the generic fiber of the numerically trivial fibration of $S$ for the canonical line bundle $K_S$ of $S$.  The handling of this problem will be the core step in this Part II.

\bigbreak\noindent{\bf\S5. Curvature of Zeroth
Direct Image of Relative Canonical and Pluricanonical Bundle.}

\bigbreak As explained at the beginning of this note the investigation of the monodromy of the zeroth direct image of the relative canonical line bundle could be traced all the way to the work of Schwarz [Schwarz1873] and Picard [Picard1885].  The differential version of the monodromy in the context of the curvature of the zeroth direct image of the relative canonical line bundle in modern mathematical terminology can be found extensively in the literature since the nineteen eighties, for example, [Viehweg1980, Kawamata1982, Viehweg-Zuo2003, Berndtsson-Paun2008, Berndtsson2006, Berndtsson2009, Schumacher2009].  Though various presentations in the existing literature seem very different, except the very different situation of the curvature of the Weil-Petersson metric from the K\"ahler-Einstein metric [Schumacher2009], eventually all the presentations come from the argument of the variation of Hodge structure about the curvature of a holomorphic vector subbundle of a flat vector bundle whose metric has both positive and negative eigenvalues.

\medbreak The key point is that the curvature of the holomorphic vector subbundle depends on whether the metric of the ambient flat vector bundle is positive or negative definite on the subspace which the second fundamental form has values in.  When the metric of the ambient flat vector bundle is defined by the integration of the exterior product of a form and its complex conjugate, the positive or negative definiteness of the metric depends on the type of the form in the subspace and on the order of primitivity of the form [Weil1958, p.22].  One important step is that very surprisingly in the case of the curvature of the zeroth direct image of the relative canonical line bundle the consideration of primitivity turns out to play no part at all.

\medbreak We present a version of the standard arguments of the second fundamental form of a subbundle of a flat bundle and the r\^ole of the type of the form in the positive or negative definiteness of the metric of integrating the exterior product of the form and its complex-conjugate.  Our purpose is to understand more precisely and explicitly, even in local coordinates, when the curvature of the zeroth direct image of the relative canonical line bundle is strictly positive.  In order to have the formulations and the notations we need, we start out with the very elementary notion of a Hodge filtration by considering a holomorphic family $\pi:X\to\Delta$ of compact complex algebraic manifolds of complex dimension $n$ over the open unit $1$-disk $\Delta$ whose fiber over the point $t\in\Delta$ we denote by $X_t$.
Of course for the introduction of the Hodge filtration it suffices to assume that the fibers are K\"ahler instead of algebraic, but for our purpose there is no advantage at all to use the more general setting and we just use the assumption of the fibers are algebraic.

\bigbreak\noindent(5.1) {\it Special Smooth Trivialization Defined by a Smooth Family of Transversal
Holomorphic Curves.}  We first introduce a special kind of smooth trivialization of $\pi:X\to\Delta$ which is
defined by a smooth family of transversal holomorphic curves.  We are only interested in the part of $X$ which is above some nonempty open neighborhood of the origin in $\Delta$.  So we will, if necessary, replace $\Delta$ by an open concentric disk of smaller positive radius.  We will do such a replacement also without explicit mention when there is no possible confusion.  We assume that $X$ is covered by a finite number of coordinate charts $U_\alpha$ with coordinates $\left(\zeta_\alpha,t\right)=\left(\zeta^1_\alpha,\cdots, \zeta^n_\alpha,t\right)$ where $t$ is the coordinate of $\Delta$.  Let $W_\alpha=X_0\cap U_\alpha$ and $\zeta_\alpha=z_\alpha$ at $t=0$ so that $X_0$ is covered by coordinate charts $W_\alpha$ with coordinates $z_\alpha=\left(z^1_\alpha,\cdots, z^n_\alpha\right)$.
By a smooth trivialization we mean a diffeomorphism
$\phi:X_0\times\Delta\to X$ so that $\phi^{-1}\left(U_\alpha\right)=W_\alpha\times\Delta$
with $W_\alpha=U_\alpha\cap X_0$.

\bigbreak\noindent(5.1.1) {\it Lemma.} After possibly replacing $\Delta$ by an open concentric disk of smaller positive radius,
we can choose a diffeomorphism
$\phi:X_0\times\Delta\to X$ such that
$\phi^*\left(z_\alpha\right)=z_\alpha\circ\phi=\zeta_\alpha\left(z_\alpha,\bar
z_\alpha,t\right)$ is holomorphic
in $t$.

\bigbreak\noindent{\it Proof.} What the lemma says is that we can
cover $X$ by a family of holomorphic curves $\left\{C_P\right\}_{P\in X_0}$
smoothly parametrized by $P\in X_0$ so
that each holomorphic curve $C_P$ is transversal to $X_0$ at $P$ and the
projection from $C_P$ to $\Delta$ is a biholomorphism (possibly after
replacing $\Delta$ by a smaller disk centered at the origin).  One
way to get this is to observe that the normal bundle of $X_0$ is
trivial and we use a global nowhere vanishing holomorphic section
$\sigma\in\Gamma\left(X_0, T_X/T_{X_0}\right)$ with $$(d\pi)(\sigma)=\left(\frac{\partial}{\partial t}\right)_{t=0}$$ to construct a smooth family of holomorphic curves $C_P$
($P\in X_0$) so that the tangent to $C_P$ at $P$ induces the value of
section $\sigma$ at $P$.  When the coordinate of $P$ in $X_0$ (with respect to the coordinate
chart $W_\alpha=U_\alpha\cap X_0$ of $X_0$)
is $z_\alpha$, the biholomorphic map from $C_P$ gives the coordinate $\left(\zeta_\alpha,t\right)$ of the point in $C_P$ above the point $t\in\Delta$
(with respect to the coordinate chart $U_\alpha$ of $X$) by the formula
$\zeta_\alpha\left(z_\alpha,\overline{z_\alpha},t\right)$. Here we choose the notation $\zeta_\alpha\left(z_\alpha,\overline{z_\alpha},t\right)$ (with the use of $\overline{z_\alpha}$ also as a variable) to describe the smooth trivialization (in terms of the coordinates $z_\alpha$ of $W_\alpha$ and the coordinates $\left(\zeta_\alpha,t\right)$ of $U_\alpha$) in order to emphasize that the function
$\zeta_\alpha\left(z_\alpha,\overline{z_\alpha},t\right)$ is only smooth in the variables $z_\alpha$ though it is holomorphic in $t$.  Later when it is too clumsy to include $\overline{z_\alpha}$ as a variable in addition to $z_\alpha$, we will drop it.  Q.E.D.

\bigbreak\noindent(5.2) {\it Kodaira-Spencer Class.}  The transformation formula of coordinates in $X$ gives
$$
\zeta_\alpha(z_\alpha,\overline{z_\alpha},t)=h_{\alpha,\beta}\left(\zeta_\beta(z_\beta,\overline{z_\beta},t),t\right)\leqno{(5.2.1)}
$$
on $U_\alpha\cap U_\beta$ with $h_{\alpha,\beta}\left(\zeta_\alpha,t\right)$ holomorphic in
both $\zeta_\alpha$ and $t$.  We differentiate (5.2.1) to get
$$
\frac{\partial\zeta^j_\alpha}{\partial t}(z_\alpha,\overline{z_\alpha},t)=\frac{\partial h^j_{\alpha,\beta}}{\partial t}\left(\zeta_\beta(z_\beta,\overline{z_\beta},t),t\right)+
\sum_{k=1}^n\frac{\partial h^j_{\alpha,\beta}}{\partial\zeta^k_\beta}\left(\zeta_\beta(z_\beta,\overline{z_\beta},t),t\right)\frac{\partial \zeta^k_\beta} {\partial t},
$$
then evaluate both sides at $t=0$ and apply $\sum_{j=1}^n\frac{\partial}{\partial z^J_\alpha}$ and then apply $\bar\partial$ to both sides to get
$$
\bar\partial\left(\sum_{j=1}^n\left(\frac{\partial\zeta^j_\alpha}{\partial t}\right)_{t=0}\frac{\partial}{\partial z^j_\alpha}\right)=\bar\partial\left(
\sum_{j,k=1}^n\left(\frac{\partial \zeta^k_\beta} {\partial t}\right)_{t=0}\frac{\partial z^j_\alpha}{\partial z^k_\beta}\frac{\partial}{\partial z^j_\alpha}\right),
$$
which is a globally defined $\bar\partial$-closed $T_{X_0}$-valued $(0,1)$-form on $X_0$.  It defines the {\it Kodaira-Spencer class}
which is an element of $H^1\left(X_0, T_{X_0}\right)$ and describes the infinitesimal deformation for the family $\pi:X\to\Delta$ at $t=0$.

\medbreak Geometrically $\sum_{j=1}^n\left(\frac{\partial\zeta^j_\alpha}{\partial t}\right)_{t=0}\frac{\partial}{\partial z^j_\alpha}$ measures on one coordinate chart $U_\alpha$ in $X$ at $t=0$ the deviation of the holomorphic coordinate $\zeta^j_\alpha$ from the background smooth trivialization coordinate $z^j_\alpha$ as a first-order infinitesimal in $t$.  Taking $\bar\partial$ to get the Kodaira-Spencer class gives the discrepancy of this measurement from one coordinate chart to another.  Note that though the Kodaira-Spencer class is obtained by taking $\bar\partial$, it is in general not zero as an element of $H^1\left(X_0, T_{X_0}\right)$, because the $\bar\partial$ is taken of a {\it local} smooth $(1,0)$-vector field even though the results of $\bar\partial$ piece together to be a global $(1,0)$-vector field valued $(0,1)$-form on $X_0$.

\medbreak In terms of $\left(z^1_\alpha,\cdots,z^n_\alpha,t\right)$ from
$$
d\zeta^j_\alpha=\sum_{k=1}^n\frac{\partial\zeta^j_\alpha}{\partial z^k_\alpha}\,dz^k_\alpha+\sum_{k=1}^n\frac{\partial\zeta^j_\alpha}{\partial\overline{z^k_\alpha}}\,d\overline{z^k_\alpha}\leqno{(5.2.2)}
$$
for any fixed $t$, we obtain the following two statements.

\medbreak\noindent(5.2.3) The $(0,1)$-component of the $(1,0)$-vector-valued $1$-form $\sum_{j=1}^n\left(\frac{\partial}{\partial t}\left(d\zeta^j_\alpha\right)\right)_{t=0}\frac{\partial}{\partial z^j_\alpha}$ equals $\bar\partial\left(\sum_{j=1}^n\left(\frac{\partial\zeta^j_\alpha}{\partial t}\right)_{t=0}\frac{\partial}{\partial z^j_\alpha}\right)$
and therefore
represents the Kodaira-Spencer class of $\pi:X\to\Delta$ at $t=0$.

\medbreak\noindent(5.2.4) $\frac{\partial}{\partial \bar t\,}\left(d\zeta^j_\alpha\right)$ vanishes, because $\zeta^j_\alpha$ is holomorphic in $t$ when expressed in terms of $\left(z^1_\alpha,\cdots,z^n_\alpha,t\right)$.

\bigbreak\noindent(5.3) {\it Variation of Hodge Structure and
the Gauss-Manin Connection.} The inclusion map $R^q\pi_*\left(\Omega^p_{X|\Delta}\right)\to R^{p+q}\pi_*\left({\mathbb C}_{X|\Delta}\right)$ of direct images over $\Delta$ is $H^q\left(X_t,\Omega^p_{X_t}\right)\to
H^{p+q}\left(X_t,{\mathbb C}\right)$ at the fiber level.  The subbundle $R^q\pi_*\left(\Omega^p_{X|\Delta}\right)$ with the induced structure from the flat vector bundle $R^{p+q}\pi_*\left({\mathbb C}_{X|\Delta}\right)$ over $\Delta$ is in general not holomorphic, but according to Griffiths [Griffiths1968]
$\bigoplus_{\ell=0}^qR^\ell\pi_*\left(\Omega^{p+q-\ell}_{X|\Delta}\right)$ is always a holomorphic subbundle of $R^{p+q}\pi_*\left({\mathbb C}_{X|\Delta}\right)$ when $0\leq\ell\leq p+q$ for the following reason.

\medbreak A $\bar\partial$-closed $(p,q)$-form $\sigma(t)$ on $X_t$ parametrized by $t$ can be written as
$$
\sigma(t)=\sum_{i_1,\cdots,i_p,\bar j_1,\cdots,\bar
j_q}\sigma^{(\alpha)}_{i_1,\cdots,i_p,\bar j_1,\cdots,\bar j_q}(t)\,
d\zeta^{i_1}_\alpha\wedge\cdots\wedge d\zeta^{i_p}_\alpha\wedge
\overline{d\zeta^{j_1}_\alpha}\wedge\cdots\wedge
\overline{d\zeta^{j_q}_\alpha}
$$
on the coordinate chart $U_\alpha$ of $X$.
Using the Gauss-Manin connection of $R^{p+q}\pi_*\left({\mathbb C}_{X|\Delta}\right)$ to differentiate with respect $t$ the class represented by $\sigma(t)$ simply means expressing $\sigma(t)$ first in terms
of $z^j_\alpha,\,\overline{z^j_\alpha},\,dz^j_\alpha,\,d\overline{z^j_\alpha}$ by (5.2.2) and then applying differentiation with respect to $t$ (for fixed $z^j_\alpha,\,\overline{z^j_\alpha}$) and finally taking the class represented by the result of differentiation.  From (5.2.3) and (5.2.4) it follows that the result of differentiating $\sigma(t)$ with respect to $\frac{\partial}{\partial\bar t\,}$ of type $(0,1)$ is a sum of a $(p,q)$-form and a $(p+1,q-1)$-form.  Hence using the Gauss-Manin connection to differentiate a local section of $\bigoplus_{\ell=0}^qR^\ell\pi_*\left(\Omega^{p+q-\ell}_{X|\Delta}\right)$ with respect to $\frac{\partial}{\partial\bar t\,}$ of type $(0,1)$ results always in a local section of $\bigoplus_{\ell=0}^qR^\ell\pi_*\left(\Omega^{p+q-\ell}_{X|\Delta}\right)$, making
$\bigoplus_{\ell=0}^qR^\ell\pi_*\left(\Omega^{p+q-\ell}_{X|\Delta}\right)$ a holomorphic subbundle of $R^{p+q}\pi_*\left({\mathbb C}_{X|\Delta}\right)$.

\bigbreak\noindent(5.4) {\it Second Fundamental Form of Zeroth Direct Image of Relative Canonical Bundle.}  To compute the second fundamental form of the zeroth direct image $R^0\pi_*\left(K_{X|\Delta}\right)$ of the relative canonical line bundle of $\pi:X\to\Delta$ as a subbundle of the flat vector bundle $R^n\pi_*\left({\mathbb C}_{X|\Delta}\right)$, we take a smooth section $\sigma(t)$ of the vector bundle associated to the locally free sheaf $R^0\pi_*\left(\Omega_{X|\Delta}^n\right)$ and for every $t$ we naturally identify $\sigma(t)$ with a holomorphic $n$-form on the fiber $X_t$.  Write $\sigma(t)=\sigma_{1,2,\cdots,n}^{(\alpha)}(t)\,
d\zeta^1_\alpha\wedge\cdots\wedge d\zeta^n_\alpha$ and use the Gauss-Manin connection to differentiate it with respect to $\frac{\partial}{\partial t\,}$ of type $(1,0)$ at $t=0$ to conclude from (5.2.3) that
$\left.\frac{\partial}{\partial t}\sigma(t)\right|_{t=0}$ is equal to the sum of some $(n,0)$-form and $$\sigma_{1,2,\cdots,n}^{(\alpha)}(0)\,\sum_{\nu=1}^n dz_\alpha^1\wedge\cdots\wedge dz_\alpha^{\nu-1}\wedge\bar\partial\left(\frac{\partial\zeta^\nu_\alpha}{\partial t}\right)_{t=0}\wedge dz^{\nu+1}_\alpha\wedge\cdots\wedge dz^n_\alpha.\leqno{(5.4.1)}
$$
For the computation of the second fundamental form we have to remove the part which is in the subbundle itself.  What is left is orthogonal to the subbundle by type considerations.  Hence the second fundamental form is given by (5.4.1), which is the contraction of the $(n,0)$-form $\sigma(0)$ with the Kodaira-Spencer class of $\pi:X\to\Delta$ at $t=0$ because of (5.2.3).  To make sure that the integral of the exterior product of the $(n-1,1)$-form in (5.4.1) and its complex conjugate gives us the correct sign after the natural normalization by a universal constant depending only on $n$, we have to check that it is in a primitive $(n-1,1)$-class in the sense that when its harmonic representative with respect to a K\"ahler form on $X_0$ is contracted with the K\"ahler form the result is zero.   This we are going to verify in the next paragraph (5.5).

\medbreak After the verification of primitivity in (5.5) the correct sign for the curvature from the second fundamental form is checked in (5.8) below.  In this note we redo the computation of our special case in order to understand more precisely and explicitly the condition for the strictness of the semipositivity of the zeroth direct image of the relative canonical line bundle and also that of the relative pluricanonical line bundle.

\bigbreak\noindent(5.5) {\it Lemma (Vanishing of Product of K\"ahler Class and Contraction of Top-Degree Holomorphic Form with Kodaira-Spencer Class).}  Let $Y$ be a compact complex K\"ahler manifold of complex dimension $n$ with K\"ahler form $\omega_0$.  Let $\varphi\in H^1\left(Y,T_Y\right)$ be a Kodaira-Spencer class of $Y$ in the sense that it represents the tangent vector of a holomorphic family of compact complex manifold over an open neighborhood of the origin in ${\mathbb C}$ whose fiber at the origin is $Y$.  Let $\sigma$ be a holomorphic $n$-form on $Y$.  Let $\psi\in H^1\left(Y,\Omega_Y^{n-1}\right)$ be defined by the contraction of $\sigma$ with $\varphi$ in the sense that if $\sigma$ is represented by $\sigma_{\alpha_1,\cdots,\alpha_n}$ in the local coordinate system $z_1,\cdots,z_n$ which is skew-symmetric in $\alpha_1,\cdots,\alpha_n$ and if $\varphi$ is represented by a $\bar\partial$-closed $(1,0)$-vector-valued $(0,1)$-form $\varphi^\alpha_{\bar\beta}$, then $\psi$ is represented by the $(n-1,1)$-form whose components with respect to the local coordinate system $z_1,\cdots,z_n$ is $\psi_{\alpha_1,\cdots,\alpha_{n-1},\bar\beta}=\sum_{\gamma=1}^n \sigma_{\alpha_1,\cdots,\alpha_{n-1},\gamma}\varphi^\gamma
_{\bar\beta}$.  Then the element of $H^2\left(Y,K_Y\right)$ which is the cup product of $\psi$ and the class represented by $\omega_0$ vanishes.  In other words, the element of $H^2\left(Y,K_Y\right)$ which is represented by the exterior product of the $(n-1,1)$-form $\psi_{\alpha_1,\cdots,\alpha_{n-1},\bar\beta}$ and the K\"ahler form $\omega_0$ vanishes.

\medbreak\noindent{\it Proof.}   Since $\varphi$ is a Kodaira-Spencer class in the sense described in the statement of the Lemma, there
exists a holomorphic family $\pi:X\to\Delta$ of compact
complex manifolds of complex dimension $n$ over the open unit $1$-disk $\Delta$ such that $Y$ is the fiber $X_0$ of the family $\pi:X\to\Delta$ at $t=0$ and $\varphi$ represents the tangent vector $\left(\frac{\partial}{\partial t}\right)_{t=0}$, where $t$ is the coordinate of $\Delta$.  By the result of Kodaira-Spencer on the stability of K\"ahler structure under holomorphic deformation [Kodaira-Spencer1960], we can assume without loss of generality that there is K\"ahler form $\omega$ on $X$ whose pullback to $Y=X_0$ is the given K\"ahler form $\omega_0$ of $Y$.

\medbreak For every tangent vector $\eta$ of type $(1,0)$ of $\Delta$ at a point $t\in\Delta$ there is a Kodaira-Spencer class $\theta\left(\eta\right)\in H^1\left(X_t,T_{X_t}\right)$ corresponding to it, where $X_t=\pi^{-1}(t)$.  The assignement $\eta\mapsto\theta(\eta)$ defines an element of $\Gamma\left(\Delta,R^1\pi_*\left(T_{X|\Delta}\right)\otimes K_\Delta\right)$ which we call the {\it
global Kodaira-Spencer class} for family $\pi:X\to\Delta$.  Note that $\varphi$ is equal to the value of $\theta$ at $\xi=\left(\frac{\partial}{\partial t}\right)_{t=0}$.

\medbreak Consider the exact sequence
$$
0\to\Omega^{p-1}_{X|\Delta}\otimes\pi^*\left(K_\Delta\right)\to\Omega^p_X\to
\Omega^p_{X|\Delta}\to 0,\leqno{(5.5.1)}
$$
which simply
says that a local holomorphic $p$-form on $X$ which pulls back to
zero on every fiber must contain as a factor the $\pi$-pullback of a
local holomorphic $1$-form on $\Delta$.  For any open subset $U$ of $\Delta$ the connecting homomorphism
$$\widehat{\Psi_{p,q-1,U}}: H^{q-1}\left(\pi^{-1}(U),\Omega^p_{X|\Delta}\right)\to
H^q\left(\pi^{-1}(U),\Omega^{p-1}_{X|\Delta}\otimes\pi^*\left(K_\Delta\right)\right)$$
from the short exact sequence $(5.5.1)$ induces the homomorphism
$$\widetilde{\Psi_{p,q-1,U}}: \Gamma\left(U, R^{q-1}\pi_*\left(\Omega^p_{X|\Delta}\right)\right)\to\Gamma\left(U,
R^q\pi_*\left(\Omega^{p-1}_{X|\Delta}\otimes\pi^*\left(K_\Delta\right)\right)\right),$$
which is the same as the homomorphism defined by contraction with the global Kodaira-Spencer class $\theta\in\Gamma\left(\Delta,R^1\pi_*\left(T_{X|\Delta}\right)\otimes K_\Delta\right)$, because

\medbreak\noindent(i) the connecting homomorphism $\widetilde{\Psi_{p,q-1,U}}$ is a measurement, as a cohomology class, of the failure of lifting an element $\tau$ of $H^{q-1}\left(\pi^{-1}(U),\Omega^p_{X|\Delta}\right)$ to an element of $H^{q-1}\left(\pi^{-1}(U),\Omega^p_X\right)$ by piecing together the local liftings of a global $\bar\partial$-closed  $\Omega^p_{X|\Delta}$-valued $(0,q-1)$-form representing $\tau$ in order to form a global $\bar\partial$-closed $\Omega^p_X$-valued $(0,q-1)$-form on $X$, and

\medbreak\noindent(ii) the Kodaira-Spencer class is precisely the
cohomology class which measures over a fiber the
discrepancy of piecing together locally holomorphic sections of $T_X$ which are the liftings of $\frac{\partial}{\partial t}$ on $\Delta$ via $\pi:X\to\Delta$.

\medbreak\noindent As the open subset $U$ of $\Delta$ varies, the map $\widetilde{\Psi_{p,q-1,U}}$ defines a sheaf-homomorphism
$$\Psi_{p,q-1}: R^{q-1}\pi_*\left(\Omega^p_{X|\Delta}\right)\to
R^q\pi_*\left(\Omega^{p-1}_{X|\Delta}\otimes\pi^*\left(K_\Delta\right)\right)$$
over $\Delta$.  After the above introduction of notations and terminology there remain only two simple steps in this proof.  The first step is the commutativity of the following diagram
$$
\begin{matrix}R^{q-1}\pi_*\left(\Omega^p_{X|\Delta}\right)&\stackrel{\Psi_{p,q-1}}{\longrightarrow}
&R^q\pi_*\left(\Omega^{p-1}_{X|\Delta}\otimes\pi^*\left(K_\Delta\right)\right)\cr
&&\cr\Phi_{p,q-1}\downarrow\qquad&&\qquad\downarrow\Phi_{p-1,q}\cr
&&\cr
R^q\pi_*\left(\Omega^{p+1}_{X|\Delta}\right)&\stackrel{\Psi_{p+1,q}}{\longrightarrow}&R^{q+1}\pi_*\left(\Omega^p_{X|\Delta}\otimes\pi^*\left(K_\Delta\right)\right)\cr\end{matrix}$$
where the map $$\Phi_{\mu,\nu}: R^\nu\pi_*\left(\Omega^\mu_{X|\Delta}\right)\to R^{\nu+1}\pi_*\left(\Omega^{\mu+1}_{X|\Delta}\right)$$
is defined by taking the cup product with the element of $H^1\left(X,\Omega^1_X\right)$ defined by the K\"ahler form $\omega$ on $X$.  This commutativity just comes from recalling the definitions of the four maps in the diagram and using the fact that the K\"ahler form $\omega$ on $X$ is $d$-closed and in particular both $\partial$-closed and $\bar\partial$-closed.

\medbreak The second step is that for the case $(p,q)=(n,1)$ the composite map
$\Psi_{n+1,1}\circ\Phi_{n,0}$ of the bundles
$$R^0\pi_*\left(\Omega^n_{X|\Delta}\right)\stackrel{\Phi_{n,0}}{\longrightarrow}R^1\pi_*\left(\Omega^{n+1}_{X|\Delta}\right)
\stackrel{\Psi_{n+1,1}}{\longrightarrow}R^2\pi_*\left(\Omega^n_{X|\Delta}\right)$$
obtained by going down the left vertical arrow first and then going to the right along the bottom horizontal arrow in the diagram is the zero map, simply because the target space $R^1\pi_*\left(\Omega^{n+1}_{X|\Delta}\right)$ of $\Phi_{n,0}$ is the zero space $\left\{0\right\}$.  Now by the commutativity of the diagram the composite map $\Phi_{n-1,1}\circ\Psi_{n,0}$ of the bundles
$$R^0\pi_*\left(\Omega^n_{X|\Delta}\right)\stackrel{\Psi_{n,0}}{\longrightarrow} R^1\pi_*\left(\Omega^{n+1}_{X|\Delta}\otimes\pi^*\left(K_\Delta\right)\right)\stackrel{\Phi_{n-1,1}}{\longrightarrow}R^2\pi_*\left(\Omega^n_{X|\Delta}\right)$$
obtained by going to the right along the top horizontal arrow first and then going down the right vertical arrow in the diagram must also be the zero map.  Now the statement of the lemma follows, because the value of ${\Psi_{n,0}}$ at the element $\sigma$ of the fiber of $R^0\pi_*\left(\Omega^n_{X|\Delta}\right)$ at $t=0$ gives us the contraction $\psi$ of the Kodaira-Spencer class $\varphi$ with $\sigma$ and the value of $\Phi_{n-1,1}$ at the contraction $\psi$ is the cup product of $\omega_0$ with $\psi$.  Q.E.D.

\bigbreak\noindent(5.6) {\it Remark on Replacing K\"ahler Class by More General (1,1)-Class}.  For the proof of Lemma (5.5) the property of $\omega_0$ being K\"ahler is not essential.  Lemma (5.5) holds if $\omega_0$ is replaced by by any $d$-closed $(1,1)$-form on $Y$ which is the pullback of some $d$-closed $(1,1)$-form on $X$ when $Y$ is a fiber of a holomorphic family $\pi:X\to\Delta$ at $t=0$.  Though Lemma (5.5) holds in such a general context, yet its known significant use so far is limited to the case of $\omega_0$ being the K\"ahler class of $Y$.

\bigbreak\noindent(5.7) {\it Equivalent Formulation of Primitivity for (n-1,1)-Forms.}  For an $(n-1,1)$-form $\psi$ with components $\psi_{\alpha_1,\cdots,\alpha_{n-1},\bar\beta}$ with respect to the local coordinate system $z_1,\cdots,z_n$ the following two formulations of its being {\it primitive} (with respect to the K\"ahler form
$\omega_0=\sum_{\lambda=1}^n dz_\lambda\wedge d\overline{z_\lambda}$) are equivalent.
One is the vanishing of the contraction of one $(1,0)$-index and one $(0,1)$-index of $\psi$ with respect to the K\"ahler form $\omega_0$.  Another is the vanishing of the exterior product $\psi\wedge\omega_0$.  The first formulation is used in [Weil1958, p.22] and the second formulation is used in Lemma (5.5).  Both agree as is shown by the following simple straightforward verification by linear algebra.   For any $1\leq\alpha_1<\cdots<\alpha_{n-2}<n$ and $\left\{\alpha_1,\cdots,\alpha_{n-1},\lambda,\nu\right\}=\left\{1,\cdots,n\right\}$ with $\lambda<\nu$ one has
$$
\sum_{\mu=1}^n\psi_{\alpha_1,\cdots,\alpha_{n-2},\mu,\bar\mu}=
(-1)^{n-\nu}\psi_{1,\cdots,\lambda-1,\lambda+1,\cdots,n,\overline{\nu}}+(-1)^{n-\lambda-1}\psi_{1,\cdots,\nu-1,\nu+1,\cdots,n,\overline{\lambda}}
$$
whereas the coefficient of $$dz_1\wedge\cdots\wedge dz_{\lambda-1}\wedge dz_{\lambda+1}\wedge\cdots\wedge dz_{\nu-1}\wedge dz_{\nu+1}\wedge\cdots\wedge dz_n\wedge\left(d\overline{z_\lambda}\wedge d\overline{z_{\nu}}\right)$$ in the $(n,2)$-form
$$
\psi\wedge\omega_0=\left(\sum_{\gamma_1,\cdots,\gamma_{n-2},\beta=1}^n\psi_{\gamma_1,\cdots,\gamma_{n-1},\bar\beta}\,dz_{\gamma_1}\wedge
\cdots\wedge dz_{\gamma_{n-1}}\wedge d\overline{z_{\beta}}\right)\wedge\left(\sum_{\mu=1}^n dz_\mu\wedge d\overline{z_\mu}\right)
$$
is equal to $(n-1)!$ times
$$
(-1)^{n-\nu}\psi_{1,\cdots,\lambda-1,\lambda+1,\cdots,n,\overline{\nu}}+(-1)^{n-\lambda-1}\psi_{1,\cdots,\nu-1,\nu+1,\cdots,n,\overline{\lambda}}.
$$
Hence the vanishing of the contraction $\sum_{\mu=1}^n\psi_{\alpha_1,\cdots,\alpha_{n-2},\mu,\bar\mu}$ is equivalent to the vanishing of $\psi\wedge\omega_0$.  The simple straightforward verification by linear algebra given here is just a very special case of the theory of exterior algebra of Hermitian spaces presented in Chapter 1 of [Weil1958].

\bigbreak\noindent(5.8) {\it Exterior Product of Primitive (n-1,1)-Form with its Complex Conjugate as Absolute-Value Square.}  As mentioned at the end of (5.4) we are going to write down explicitly the simple straightforward linear algebra to understand more precisely the condition for the strictness of the semipositivity of the zeroth direct image of the relative canonical line bundle and also that of the relative pluricanonical line bundle.  For a primitive $(n-1,1)$-form $\psi$ given by
$$\psi=\sum_{\alpha_1,\cdots,\alpha_{n-1},\beta=1}^n \psi_{\alpha_1,\cdots,\alpha_{n-1}\bar\beta}\,dz_{\alpha_1}\wedge\cdots\wedge dz_{\alpha_{n-1}}\wedge d\overline{z_\beta}$$
the expression $\frac{1}{\left((n-1)!\right)^2}\,\psi\wedge\overline{\psi}$ is equal to
$\sum_{\lambda,\nu=1}^n\left|\psi_{1,\cdots,\lambda-1,\lambda+1,\cdots,n,\overline{\nu}}\right|^2$
times
$\left(dz_1\wedge\cdots\wedge dz_n\right)\wedge\overline{\left(dz_1\wedge\cdots\wedge dz_n\right)}$.

\medbreak This is verified by showing that in general, without assuming that the $(n-2,2)$-form $\psi$ is primitive, the expression $\frac{1}{\left((n-1)!\right)^2}\,\psi\wedge\overline{\psi}$ is equal to
$$
-\frac{1}{\left((n-2)!\right)^2}\sum_{1\leq\alpha_1<\cdots<\alpha_{n-2}\leq n}\left|\sum_{\mu=1}^n \psi_{\alpha_1,\cdots,\alpha_{n-2},\mu,\bar\mu}\right|^2+\sum_{\lambda,\nu=1}^n\left|\psi_{1,\cdots,\lambda-1,\lambda+1,\cdots,n,\overline{\nu}}\right|^2$$
times
$\left(dz_1\wedge\cdots\wedge dz_n\right)\wedge\overline{\left(dz_1\wedge\cdots\wedge dz_n\right)}$.

\medbreak Here is the verification.  The expression $\psi\wedge\overline{\psi}$ is equal to
$$
\sum_{1\leq\alpha_1,\cdots,\alpha_{n-1},\beta\leq n\atop
1\leq\gamma_1,\cdots,\gamma_{n-1},\delta\leq n}\psi_{\alpha_1,\cdots,\alpha_{n-1}\bar\beta}\,dz_{\alpha_1}\wedge\cdots\wedge dz_{\alpha_{n-1}}\wedge d\overline{z_\beta}\wedge\overline{
\left(\psi_{\gamma_1,\cdots,\gamma_{n-1}\bar\delta}\,dz_{\gamma_1}\wedge\cdots\wedge dz_{\gamma_{n-1}}\wedge d\overline{z_\delta}\right)}.
$$
For a permutation $\alpha_1,\cdots,\alpha_n$ of $1,\cdots,n$ let $\epsilon_{\alpha_1,\cdots,\alpha_n}$ be its signature.  For a term
$$\psi_{\alpha_1,\cdots,\alpha_{n-1}\bar\beta}\,dz_{\alpha_1}\wedge\cdots\wedge dz_{\alpha_{n-1}}\wedge d\overline{z_\beta}\wedge\overline{
\left(\psi_{\gamma_1,\cdots,\gamma_{n-1}\bar\delta}\,dz_{\gamma_1}\wedge\cdots\wedge dz_{\gamma_{n-1}}\wedge d\overline{z_\delta}\right)}$$
to be nonzero, we must have the re-arrangements $\alpha_1,\cdots,\alpha_n$ and $\gamma_1,\cdots,\gamma_n$ of $1,\cdots,n$ and must have $\beta=\gamma_n$ and $\delta=\alpha_n$ so that we have
$$
\displaylines{\psi_{\alpha_1,\cdots,\alpha_{n-1}\bar\beta}\,dz_{\alpha_1}\wedge\cdots\wedge dz_{\alpha_{n-1}}\wedge d\overline{z_\beta}\wedge\overline{
\left(\psi_{\gamma_1,\cdots,\gamma_{n-1}\bar\delta}\,dz_{\gamma_1}\wedge\cdots\wedge dz_{\gamma_{n-1}}\wedge d\overline{z_\delta}\right)}
\cr=\psi_{\alpha_1,\cdots,\alpha_{n-1}\overline{\gamma_n}}dz_{\alpha_1}\wedge\cdots\wedge dz_{\alpha_{n-1}}\wedge d\overline{z_{\gamma_n}}\wedge\overline{
\left(\psi_{\gamma_1,\cdots,\gamma_{n-1}\overline{\alpha_n}}dz_{\gamma_1}\wedge\cdots\wedge dz_{\gamma_{n-1}}\wedge d\overline{z_{\alpha_n}}\right)}\cr
\qquad=-\epsilon_{\alpha_1,\cdots,\alpha_n}\,\epsilon_{\gamma_1,\cdots,\gamma_n}\,
\psi_{\alpha_1,\cdots,\alpha_{n-1}\overline{\gamma_n}}\hfill\cr\hfill\overline{
\psi_{\gamma_1,\cdots,\gamma_{n-1}\overline{\alpha_n}}}\left(dz_1\wedge\cdots\wedge dz_n\right)\wedge\overline{\left(dz_1\wedge\cdots\wedge dz_n\right)}.\qquad}
$$
For our computation we need only consider the case where $\alpha_1,\cdots,\alpha_{n-1}$ and $\gamma_1,\cdots,\gamma_{n-1}$ are strictly increasing and $\alpha_n=\lambda$ and $\gamma_n=\nu$.  Then we get
$$
\displaylines{\psi_{\alpha_1,\cdots,\alpha_{n-1}\bar\beta}\,dz_{\alpha_1}\wedge\cdots\wedge dz_{\alpha_{n-1}}\wedge d\overline{z_\beta}\wedge\overline{
\left(\psi_{\gamma_1,\cdots,\gamma_{n-1}\bar\delta}\,dz_{\gamma_1}\wedge\cdots\wedge dz_{\gamma_{n-1}}\wedge d\overline{z_\delta}\right)}
\cr=-(-1)^{\lambda+\nu}\psi_{1,\cdots,\lambda-1,\lambda+1,\cdots,n,\overline{\nu}}\,
\overline{\psi_{1,\cdots,\nu-1,\nu+1,\cdots,n,\overline{\lambda}}}
\left(dz_1\wedge\cdots\wedge dz_n\right)\wedge\overline{\left(dz_1\wedge\cdots\wedge dz_n\right)}.}
$$
Then $\frac{1}{\left((n-1)!\right)^2}\,\psi\wedge\overline{\psi}$ is equal to
$$
\sum_{\lambda,\nu=1}^n(-1)^{\lambda+\nu}\psi_{1,\cdots,\lambda-1,\lambda+1,\cdots,n,\overline{\nu}}\,
\overline{\psi_{1,\cdots,\nu-1,\nu+1,\cdots,n,\overline{\lambda}}}
\left(dz_1\wedge\cdots\wedge dz_n\right)\wedge\overline{\left(dz_1\wedge\cdots\wedge dz_n\right)}.
$$
For any $1\leq\alpha_1<\cdots<\alpha_{n-2}<n$ we let $\left\{\alpha_1,\cdots,\alpha_{n-1},\lambda,\nu\right\}=\left\{1,\cdots,n\right\}$ with $\lambda<\nu$ and get
$$
\sum_{\mu=1}^n \psi_{\alpha_1,\cdots,\alpha_{n-2},\mu,\bar\mu}=
(-1)^{n-\nu}\psi_{1,\cdots,\lambda-1,\lambda+1,\cdots,n,\overline{\nu}}+(-1)^{n-\lambda-1}\psi_{1,\cdots,\nu-1,\nu+1,\cdots,n,\overline{\lambda}}.
$$
Its absolute-value square gives
$$
\displaylines{\quad\left|\sum_{\mu=1}^n \psi_{\alpha_1,\cdots,\alpha_{n-2},\mu,\bar\mu}\right|^2=
\left|\psi_{1,\cdots,\lambda-1,\lambda+1,\cdots,n,\overline{\nu}}\right|^2+\left|\psi_{1,\cdots,\nu-1,\nu+1,\cdots,n,\overline{\lambda}}\right|^2\hfill\cr\hfill
+(-1)^{\lambda+\nu+1}\left(\psi_{1,\cdots,\lambda-1,\lambda+1,\cdots,n,\overline{\nu}}\,\overline{\psi_{1,\cdots,\nu-1,\nu+1,\cdots,n,\overline{\lambda}}}
+\overline{\psi_{1,\cdots,\lambda-1,\lambda+1,\cdots,n,\overline{\nu}}}\,\psi_{1,\cdots,\nu-1,\nu+1,\cdots,n,\overline{\lambda}}\right).\cr}
$$
The final conclusion is that $\frac{1}{\left((n-1)!\right)^2}\,\psi\wedge\overline{\psi}$ is equal to
$$
-\frac{1}{\left((n-2)!\right)^2}\sum_{1\leq\alpha_1<\cdots<\alpha_{n-2}\leq n}\left|\sum_{\mu=1}^n \psi_{\alpha_1,\cdots,\alpha_{n-2},\mu,\bar\mu}\right|^2+\sum_{\lambda,\nu=1}^n\left|\psi_{1,\cdots,\lambda-1,\lambda+1,\cdots,n,\overline{\nu}}\right|^2$$
times
$\left(dz_1\wedge\cdots\wedge dz_n\right)\wedge\overline{\left(dz_1\wedge\cdots\wedge dz_n\right)}$.  When the $(n-2,2)$-form $\psi$ is primitive, we have
$$\frac{1}{\left((n-1)!\right)^2}\,\psi\wedge\overline{\psi}=\left(\sum_{\lambda,\nu=1}^n\left|\psi_{1,\cdots,\lambda-1,\lambda+1,\cdots,n,\overline{\nu}}\right|^2\right)
\left(dz_1\wedge\cdots\wedge dz_n\right)\wedge\overline{\left(dz_1\wedge\cdots\wedge dz_n\right)},$$
which gives us not only the semipositivity for the curvature we need but, more importantly to us, also the explicit condition for the curvature to vanish.

\bigbreak\noindent(5.9) {\it Lemma on Strict Positivity of Zeroth Direct Image of Relative Canonical Line Bundle.}  Let $\pi:X\to\Delta$ be a holomorphic family of compact compact manifolds of complex dimension $n$ over the open unit $1$-disk $\Delta$.  Let the metric of the zeroth direct image $R^0\pi_*\left(K_{X|\Delta}\right)$ be defined by the integral of the absolute-value square of a holomorphic $n$-form on a fiber.  Let $\hat\sigma$ be the nonzero element of the fiber of $R^0\pi_*\left(K_{X|\Delta}\right)$ at the origin $t=0$ of $\Delta$ defined by a non identically zero holomorphic $n$-form $\sigma$ on the fiber $X_0$ of $\pi:X\to\Delta$ at $t=0$.  Then the following statements hold.

\medbreak\noindent(a) The curvature of the zeroth direct image $R^0\pi_*\left(K_{X|\Delta}\right)$ with respect to the metric described above is seimpositive.

\medbreak\noindent(b)  The curvature of the zeroth direct image $R^0\pi_*\left(K_{X|\Delta}\right)$ at $\hat\sigma$ is zero if and only if the element of $H^1\left(X_0, \Omega_{X_0}^{n-1}\right)$ defined by the contraction of $\sigma$ with the Kodaira-Spencer class of $\pi:X\to\Delta$ at $t=0$ is zero.

\medbreak\noindent(c)  The curvature of the zeroth direct image $R^0\pi_*\left(K_{X|\Delta}\right)$ at $\hat\sigma$ is zero if and only if the Kodaira-Spencer class of $\pi:X\to\Delta$ at $t=0$ is zero.

\bigbreak\noindent(5.9.1) {\it Proof of (5.9).}  The statements (a) and (b) follow immediately from the arguments presented from (5.4) to (5.8) above.  We would like to draw attention to one point.  The vanishing of the cup product of the K\"ahler class and the contraction of the olomorphic $n$-form with the Kodaira-Spencer class given in Lemma (5.5) is only as an element of the cohomology group.  So to get its vanishing we have to use a harmonic representative (with respect some K\"ahler form $\omega_0$ on $X_0$).  Note that both the process of contraction with the K\"ahler form $\omega_0$ and the process of taking the exterior product with $\omega_0$ commute with the process of taking harmonic components with respect to $\omega_0$.  Hence for the statement (b) the curvature is nonzero only when the harmonic representative of the contraction of $\sigma$ with the Kodaira-Spencer class is nonzero and not just when the contraction of $\sigma$ with the $(1,0)$-vector-valued $(0,1)$-form in (5.2) representing the Kodaira-Spencer class is not identically zero.

\medbreak We now turn our attention to the statement (c).  Suppose the Kodaira-Spencer class at $t=0$ for $\pi:X\to\Delta$ is nonzero.  Let it be represented by the $(1,0)$-vector-valued $\bar\partial$-closed $(0,1)$-form $\varphi$, for example, constructed in (5.2), whose components with respect to a local coordinate system of $X_0$ are $\varphi^\alpha_{\bar\beta}$.  Assume that the element in $H^1\left(X_0,\Omega_{X_0}^{n-1}\right)$ represented by the contraction of $\sigma$ and $\varphi$ is zero and we are going to derive the contradiction that $\varphi$ is equal to $\bar\partial$ of some global $(1,0)$-vector field on $X_0$.

\medbreak Fix a local coordinate system of $X_0$. Let $\psi$ be the $\bar\partial$-closed $(n-1,1)$-form which is the contraction of $\sigma$ and $\varphi$ and whose components with respect to the local coordinate system of $X_0$ are $\psi_{\alpha_1,\cdots,\alpha_{n-1},\bar\beta}$ which are skew-symmetric in $\alpha_1,\cdots,\alpha_{n-1}$.
By assumption $\psi$ is equal to $\bar\partial\theta$ for some $(n-1,0)$-form $\theta$ whose components with respect to the local coordinate system of $X_0$ are $\theta_{\alpha_1,\cdots,\alpha_{n-1}}$ which are skew-symmetric in $\alpha_1,\cdots,\alpha_{n-1}$.  Hence
$\psi_{1,\cdots,\lambda-1,\lambda+1,\cdots,n,\bar\nu}=\partial_{\bar\nu}\theta_{1,\cdots,\lambda-1,\lambda+1,\cdots,n}$.  Let the components of $\sigma$ with respect to the local coordinate system of $X_0$ be $\sigma_{\alpha_1,\cdots,\alpha_n}$ which are skew-symmetric in $\alpha_1,\cdots,\alpha_n$.
We know that
$$\displaylines{\psi_{1,\cdots,\lambda-1,\lambda+1,\cdots,n,\bar\nu}=\sum_{\gamma=1}^n\sigma_{1,\cdots,\lambda-1,\lambda+1,\cdots,\gamma-1,\gamma+1,\cdots,n,\gamma}\varphi^\gamma_{\bar\nu}\cr
=\sigma_{1,\cdots,n}\left(\sum_{1\leq\gamma<\lambda}(-1)^{\gamma-n-1}\varphi^\gamma_{\bar\nu}+\sum_{\lambda<\gamma\leq n}(-1)^{\gamma-n}\varphi^\gamma_{\bar\nu}\right).}
$$
Define $\tau^\gamma$ by
$$
\sum_{1\leq\gamma<\lambda}(-1)^{\gamma-n-1}\tau^\gamma+\sum_{\lambda<\gamma\leq n}(-1)^{\gamma-n}\tau^\gamma
=\frac{\theta_{1,\cdots,\lambda-1,\lambda+1,\cdots,n}}{\sigma_{1,\cdots,n}}.
$$
Then $\tau^\gamma$ gives a well-defined $(1,0)$-vector field whose coefficients are distributions, because the local holomorphic function $\sigma_{1,\cdots,n}$ vanishes at the zero-set of the holomorphic $n$-form $\sigma$ which may be nonempty.  Because of the possible vanishing of $\sigma_{1,\cdots,n}$, instead of the equation $\partial_{\bar\beta}\tau^\gamma=\varphi^\gamma_{\bar\beta}$ we can get in general only the weaker conclusion that $\partial_{\bar\beta}\tau^\gamma=\varphi^\gamma_{\bar\beta}+V$ for some $(1,0)$-vector-valued $(0,1)$-current $V$ whose nonsmooth property comes locally from $\bar\partial\left(\frac{1}{h}\right)$ of some local holomorphic function $h$ which is not nowhere zero.  The $(1,0)$-vector-valued $(0,1)$-current $V$ must be nonzero, otherwise we would have a contradiction from the equation $\partial_{\bar\beta}\tau^\gamma=\varphi^\gamma_{\bar\beta}$ which implies that the Kodaira-Spencer class as an element in $H^1\left(X_0,\Omega_{X_0}\right)$ defined by $\varphi^\gamma_{\bar\beta}$ vanishes.

\medbreak We observe that the multiplication of $V$ and the characteristic function of any subvariety of complex codimension $\geq 2$ in $X_0$ is zero, because the kind of singularity from $\bar\partial\left(\frac{1}{h}\right)$ cannot be supported on such a thin set as a subvariety of complex codimension $\geq 2$ in $X_0$.  At a regular point of the zero-set of $h$ where the zero-set of $h$ is given locally by the vanishing of the coordinate $z_n$ of a local coordinate system $z_1,\cdots,z_n$ of $X_0$, the value of $\bar\partial\left(\frac{1}{h}\right)$ at a compactly supported smooth $(2n-1)$-form $\rho$ is zero unless $\rho$ is of the form $\hat\rho\wedge dz_n$ for some $(n-1,n-1)$-form $\hat\rho$ in which case the value is equal to the integral of $\hat\rho$ on the local complex hypersurface defined by the vanishing of $h$ with multiplicity counted.  In other words, a non-identically-zero nonsmooth tangent-vector-valued $(0,1)$-form with nonsmoothness coming from $\bar\partial\left(\frac{1}{h}\right)$ has to come from the one obtained by raising the $(1,0)$-index of the $(1,1)$-current $\partial\bar\partial\log|h|^2$, but then the result of the index-raising would force the tangent-vector value to be of type $(0,1)$.  This shows that the $(1,0)$-vector-valued $V$ must be zero and we get a contradiction. Q.E.D.

\bigbreak\noindent(5.9.2) {\it Remark.}  Though we derive the very strong conclusion of (5.9)(c), yet in our proof of the abundance conjecture the weaker conclusion of (5.9)(b) suffices.  The conclusion of (5.9)(b) or (5.9)(c) actually will be applied to the case of the zeroth direct image of the relative pluricanonical line bundle for the case of the numerically trivial fibration (assumed to coincide with the numerically trivial foliation) for the canonical line bundle by using the technique of taking roots.  We are going to present this in the next section.

\bigbreak\noindent(5.10) {\it Strict Positivity of Zeroth Direct Image of Relative Pluricanonical Line Bundle.}  Given a holomorphic family $\pi:X\to\Delta$ of compact complex algebraic manifolds of complex dimension $n$ over the open unit $1$-disk $\Delta$, the technique of taking the $m$-root of a pluricanonical section, similar to the procedure described in (3.1.1), can be applied to the family $\pi:X\to\Delta$ outside a discrete subset of $\Delta$ to reduce the question of the curvature of the zeroth direct image $R^0\pi_*\left(mK_{X|\Delta}\right)$ of the relative $m$-canonical line bundle at a generic point of the base $\Delta$ to the question of the curvature of the zeroth direct image $R^0\tilde\pi_*\left(K_{\tilde X|\Delta}\right)$ of the relative canonical line bundle of another holomorphic family $\tilde\pi:\tilde X\to\Delta$ at a generic point of $\Delta$.  After the application of this technique of taking roots, the conclusions of Lemma (5.5) still hold at a generic point of $\Delta$ when the relative canonical line bundle $K_{X|\Delta}$ is replaced the relative $m$-canonical line bundle $mK_{X|\Delta}$.  The reason why we have to look at a generic point of $\Delta$ is that the newly constructed holomorphic family of compact complex algebraic manifolds by taking roots is a genuine regular family only over $\Delta$ minus a discrete subset of $\Delta$.

\bigbreak\noindent{\bf\S6. Strict Positivity of Direct Image of Relative Pluricanonical Bundle Along Numerically Trivial Fibers in the Base of Numerically Trivial Fibration.}

\bigbreak To finish the proof of the abundance conjecture under the assumption of the coincidence always of numerically trivial foliation and foliation for the canonical line bundle, we need to verify that, in the numerically trivial fibration for the canonical line bundle, the curvature of the direct image of the relative pluricanonical line bundle is strictly positive along a generic numerically trivial fiber in the base for the canonical line bundle of the base.  We are doing this argument here in two ways, one using the pluricanonical version of  (5.5)(b) and the other using the pluricanonical version of the stronger (5.5)(c).  One important tool which we use is the following technique of Siegel to bound the dimension of the space of holomorphic sections by Schwarz's lemma.

\bigbreak\noindent(6.1) {\it Effect of Fixed Ample Twisting on Growth of Dimension of Sections of Line Bundles.} The following proposition adapted from the technique of Carl Siegel [Siegel1958] facilitates the control of the growth of the space of sections as $m\to\infty$ when the twisting by a fixed ample line bundle is added to the $m$-canonical line bundle and when one has the vanishing of the Kodaira-Spencer class at a generic point of a complex curve in the parameter space of the family from the pluricanonical version of (5.5)(c) or when one has the pluricanonical version of the weaker (5.5)(b).  The key point is that the uniform bounds of the transition functions can be controlled by the curvature via the integration along a short curve-segment to obtain parallel sections along the curve segment and via the $L^2$ estimates of $\bar\partial$ and also controlled by the representatives of Kodaira-Spencer classes and in particular their vanishing.  The addition of a fixed ample line bundle affects the uniform bounds of the transition functions of a sequence of line bundles only by fixed amounts independent of the line bundles in the sequence.  We need the result from Siegel's technique only for the case of uniformly bounded transition functions or uniformly bounded curvature.  A general form of the result deals with the case of growth order for transition functions or curvature.  We will give the general case in (6.1.1) and (6.1.2) and then the case of uniform bound which we need in (6.1.3).

\bigbreak\noindent(6.1.1) {\it Proposition (Growth of Dimension of Space of Sections).}  Let $Y$ be a compact complex manifold of complex dimension $n$ and $\left\{U_j\right\}_{1\leq j\leq J}$ be a finite cover of $Y$ by coordinate balls.  Let $L_m$ for $m\in{\mathbb N}$ be a holomorphic line bundle with transition functions $g_{m,jk}$ on $U_j\cap U_k$ for the transformation of fiber coordinates from $U_k$ to $U_j$.  Let $\gamma<\delta$ and $C$ be positive numbers such that $\sup_{U_j\cap U_k}\left|g_{m,jk}\right|\leq C^{m^{\frac{\gamma}{n}}}$.  Then
$$
\lim_{m\to\infty}\frac{1}{m^\delta}\dim_{\mathbb C}\Gamma\left(Y,L_m\right)=0.
$$

\medbreak\noindent{\it Proof.} Assume the contrary so that there exist a sequence $m_\nu$ in $N$ and some $a>0$ such that
$$
\dim_{\mathbb C}\Gamma\left(Y,L_{m_\nu}\right)\geq a\left(m_\nu\right)^\delta.\leqno{(6.1.1.1)}
$$
We are going to derive a contradiction.  Without loss of generality we assume that the radius of each coordinate ball $U_j$ is $1$.  There exists some $0<r<1$ such that the coordinate balls $W_j$ concentric with $U_j$ with radius $r$ (for $1\leq j\leq J$) still cover $Y$.  Choose $\gamma<\eta<\delta$.  Let $P_j$ be the center of the two concentric coordinate balls $W_j$ and $U_j$.  It follows from $(6.1.1.1)$ and $\eta<\delta$ that there exists $\nu_0\in{\mathbb N}$ such that for $\nu\geq \nu_0$ we can find a nonzero element $s_{m_\nu}\in\Gamma\left(Y,L_{m_\nu}\right)$ represented by a holomorphic function $s_{m_\nu,j}$ on $U_j$ such that $s_{m_\nu,j}=g_{m_\nu,jk}s_{m_\nu,k}$ on $U_j\cap U_k$ and $s_{m_\nu,j}$ vanishes to order $\geq\left(m_\nu\right)^{\frac{\eta}{n}}$ at each of the centers $P_j$, because the number of terms of degree $\leq q$ in a power series in $n$ variables is $n+q\choose n$ which is of the order $q^n$ as $q\to\infty$.  By multiplying each $s_{m_\nu}$ by a positive constant, we can assume without loss of generality that
$$\max_{1\leq j\leq J}\sup_{U_j}\left|s_{m_\nu,j}\right|=1.
$$
Let
$$
A_\nu=\max_{1\leq j\leq J}\sup_{W_j}\left|s_{m_\nu,j}\right|.
$$
Because $s_{m_\nu,j}$ vanishes to order $\geq\left(m_\nu\right)^{\frac{\eta}{n}}$, by Schwarz's lemma $A_\nu\leq r^{\left(m_\nu\right)^{\frac{\eta}{n}}}$.  On the other hand, since each $U_j$ is contained in the union of $W_k$ for $1\leq k\leq J$ and since $\sup_{U_j\cap U_k}\left|g_{m,jk}\right|\leq C^{m^{\frac{\gamma}{n}}}$, it follows from $s_{m_\nu,j}=g_{m_\nu,jk}s_{m_\nu,k}$ on $U_j\cap U_k$ that $1\leq C^{\left(m_\nu\right)^{\frac{\gamma}{n}}}A_\nu$.  As $\nu\to\infty$ we arrive at a contradiction from
$$
1\leq C^{\left(m_\nu\right)^{\frac{\gamma}{n}}}A_\nu\leq C^{\left(m_\nu\right)^{\frac{\gamma}{n}}}r^{\left(m_\nu\right)^{\frac{\eta}{n}}}
$$
and $\eta<\gamma$, because $r<1$.  Q.E.D.

\bigbreak\noindent(6.1.2) {\it Remark on Vector Bundle Version and Version with Curvature Condition.} There is a vector bundle version of (6.1.1) which holds for a sequence of vector bundles of fixed rank over $Y$ and is obtained by considering the projectivization of the vector bundles to reduce the vector bundle case to the line bundle case.  Note that when the projectivization procedure is used, instead of one single manifold $Y$, there is a sequence of manifolds, but the estimates can be done uniformly from the data on each manifold.  There is also a version of (6.1.1) with curvature condition similar to (6.1.3)(c) below.

\bigbreak\noindent(6.1.3) {\it Proposition (Dimension of Section Space for Bundle of Bounded Curvature).}  Let $Y$ be a compact complex manifold of complex dimension $n$ and $\left\{U_j\right\}_{1\leq j\leq J}$ be a finite cover of $Y$ by coordinate balls.  Let $g$ be a Hermitian metric of $Y$.

\medbreak\noindent(a) Let $L_m$ for $m\in{\mathbb N}$ be a holomorphic line bundle with transition functions $g_{m,jk}$ on $U_j\cap U_k$ for the transformation of fiber coordinates from $U_k$ to $U_j$.  If $\sup_{y\in U_j\cap U_k, j\not=k}\left|g_{m,jk}(y)\right|<\infty$, then $\sup_{m\in{\mathbb N}}\dim_{\mathbb C}\Gamma\left(Y,L_m\right)<\infty$.

\medbreak\noindent(b)
Let $L_m$ be a sequence of holomorphic line bundles with smooth Hermitian metric $h_m$ for $m\in{\mathbb N}$.  If the curvature form of $h_m$ is uniformly bounded with respect to $g$ for $m\in{\mathbb N}$, then
$\sup_{m\in{\mathbb N}}\dim_{\mathbb C}\Gamma\left(Y,L_m\right)<\infty$.

\medbreak\noindent(c)  Let $r\in{\mathbb N}$ and $V_m$ be a sequence of holomorphic vector bundles of rank $r$ with smooth Hermitian metric $H_m$ for $m\in{\mathbb N}$.  If the curvature form of $H_m$ is uniformly bounded with respect to $g$ and $H_m$ for $m\in{\mathbb N}$, then $\sup_{m\in{\mathbb N}}\dim_{\mathbb C}\Gamma\left(Y,V_m\right)<\infty$.

\medbreak\noindent{\it Proof.} The proof of (a) is completely analogous to and actually easier than the proof of (6.1.1).  The proof of (b) uses the curvature estimate to get local trivializations of $L_m$ whose transition functions satisfy the uniform bound estimates for some finite cover by coordinate balls independent of $m$.  The proof of (c) uses the projectivization of $V_m$ to reduce the problem to a sequence of holomorphic line bundles over different complex manifolds $Y_m$ of complex dimension $n+r-1$, but the arguments for the estimates can be carried out uniformly with respect to the complex manifolds $Y_m$.

\bigbreak\noindent(6.2) {\it Amply Twisted Pluricanonical Section as Product Pluricanonical Section and Ample Bundle Section for Manifold of Zero Numerical Kodaira Dimension.}  Let $Y$ be a compact complex algebraic manifold whose numerical Kodaira dimension is zero and hence whose Kodaira dimension is also zero.  We know that for some sufficiently divisible positive integer $m_0$ we have $\dim_{\mathbb C}\Gamma\left(Y, pm_0K_Y\right)=1$ for $p\in{\mathbb N}$.  Let $s_{m_0}$ be a nonzero element of $\Gamma\left(Y, m_0K_Y\right)$.  Let $B$ be an ample line bundle over $Y$.  Then we know that $\sup_{m\in{\mathbb N}}\Gamma\left(Y, mK_Y+B\right)<\infty$.  Hence for $q_0\in{\mathbb C}$ sufficiently large $\dim_{\mathbb C}\Gamma\left(Y, pq_0m_0K_Y+B\right)$ is independent of $p\in{\mathbb N}$.  As a consequence the map from $\Gamma\left(Y, pq_0m_0K_Y+B\right)$ to $\Gamma\left(Y, (p+1)q_0m_0K_Y+B\right)$ defined by multiplication by $\left(s_{m_0}\right)^{q_0}$ is an isomorphism.  In other words, we have the following statement.

\medbreak\noindent(6.2.1) Every element of $\Gamma\left(Y, (p+1)q_0m_0K_Y+B\right)$ can be decomposed into a product of an element of $\Gamma\left(Y, pq_0m_0K_Y+B\right)$ and the global $q_0m_0$-canonical section $\left(s_{m_0}\right)^{q_0}$ of $Y$.

\medbreak\noindent(6.3) {\it Factoring of Direct Image of Amply Twisted Relative Pluricanonical Bundle in Numerically Trivial Fibration.} Consider now the numerically trivial fibration $\pi:X\to S$ for the canonical line bundle of $X$ and an ample line bundle $A$ on $X$.  By desingularization we can without loss of generality assume that $S$ is a manifold.  By (6.2.1) we can decompose the zeroth direct image $R^0\pi_*\left(pq_0m_0K_{X|S}+A\right)$ as the product of the line bundle $\left(R^0\pi_*\left(q_0m_0K_{X|S}\right)\right)^{\otimes(p-1)}$ and the vector bundle $R^0\pi_*\left(q_0m_0K_{X|S}+A\right)$.  We consider this decomposition first at a generic point of $S$ and then use the deformational invariance of twisted plurigenera to consider some other points of $S$.  Here we have left out the details of the discussion of singular fibers for $\pi:X\to S$ and, in our sketch of the proof, assume that we have the simpler case of $R^0\pi_*\left(pq_0m_0K_{X|S}+A\right)$, $\left(R^0\pi_*\left(q_0m_0K_{X|S}\right)\right)^{\otimes(p-1)}$, and $R^0\pi_*\left(q_0m_0K_{X|S}+A\right)$ all being vector bundles instead of being general coherent sheaves.

\bigbreak\noindent(6.4) {\it Strict Positivity of Direct Image of Relative Pluricanonical Bundle Along Generic Fibers of Numerical Trivial Fibration of Base.}  Let $\pi:X\to S$ be the numerically trivial fibration for $K_X$ and and $\tau:S\to T$ be the numerically trivial fibration for $K_S$.  Assume that the dimension of a generic fiber of $\tau:S\to T$ is positive ({\it i.e.,} $S$ is not of general type).  Let $A$ be an ample line bundle on $X$.  Let $k$ be any positive integer.  By using blow-ups, we can assume without loss of generality that both $\pi$ and $\tau$ are holomorphic maps and $S$ and $T$ are manifolds.  We are going to verify that for some $m\in{\mathbb N}$ the curvature of the restriction of the line bundle $R^0\pi_*\left(mK_{X|S}\right)$ to $\tau^{-1}\left(t_0\right)$ is strictly positive at a generic point $t_0$ of $T$.  Suppose the contrary.  If we use the pluricanonical version of (5.5)(b) instead of the pluricanonical version of (5.5)(c), by (6.3) we can write
$R^0\pi_*\left(pq_0m_0K_{X|S}+kA\right)$ as the tensor product of $\left(R^0\pi_*\left(q_0m_0K_{X|S}\right)\right)^{\otimes(p-1)}$, and $R^0\pi_*\left(q_0m_0K_{X|S}+kA\right)$ when we restrict all three to $\tau^{-1}\left(t_0\right)$.  Since the restriction of $\left(R^0\pi_*\left(q_0m_0K_{X|S}\right)\right)^{\otimes(p-1)}$ to $\tau^{-1}\left(t_0\right)$ is a line bundle whose curvature vanishes, it follows that the curvature of the restriction of $R^0\pi_*\left(pq_0m_0K_{X|S}+kA\right)$ to $\tau^{-1}\left(t_0\right)$ is uniformly bounded independent of $p$.   By (6.1.2) the dimension of global holomorphic sections of $R^0\pi_*\left(pq_0m_0K_{X|S}+A\right)$ over
$\tau^{-1}\left(t_0\right)$ is uniformly bounded independent of $p$.  (If we use the pluricanonical version of (5.5)(c), then this uniform bound of curvature is immediate without the use of (6.1.2).)  Since $\tau^{-1}\left(t_0\right)$ is a generic fiber of the numerically trivial fibration of $S$ for $K_S$, we conclude that the dimension of $\Gamma\left(\left(\tau\circ\pi\right)^{-1}\left(t_0\right), pq_0K_X+kA\right)$ is uniformly bounded in $p$.  It follows that $K_X$ is numerically trivial on $\left(\tau\circ\pi\right)^{-1}$, which contradicts the fact that $\pi:X\to S$ is the numerically trivial fibration of $X$ for $K_X$, because the dimension of $\left(\tau\circ\pi\right)^{-1}$ is strictly greater than the dimension of a generic fiber of $\pi:X\to S$.  Here we have suppressed some details, because a precise argument needs the finite upper bound of the dimension of $\Gamma\left(\left(\tau\circ\pi\right)^{-1}\left(t_0\right), mK_X+kA\right)$ for $m\in{\mathbb N}$ instead of just the finite upper bound of the dimension of $\Gamma\left(\left(\tau\circ\pi\right)^{-1}\left(t_0\right), pq_0K_X+kA\right)$ on $p\in{\mathbb N}$.

\bigbreak\noindent(6.5) {\it Construction of Pluricanonical Sections by $L^2$ Estimates of $\bar\partial$.}  Let $X^{(0)}=X$ and we construct inductively the numerically trivial fibration $\pi^{(\nu)}:X^{(\nu)}\to X^{(\nu+1)}$ for $K_{X^{(\nu)}}$ for $0\leq\nu<J$ so that the generic fiber dimension for $\pi^{(\nu)}: X^{(\nu)}\to X^{(\nu+1)}$ is positive and $X^{(J)}$ is either a single point or of general type.  Here we assume without loss of generality, after using blow-ups if necessary, that each $X^{(\nu)}$ is a complex algebraic manifold for $0\leq\nu\leq J$.  Also we assume that we know that numerically trivial foliation for the canonical line bundle coincides with the numerically trivial foliation for the canonical line bundle, which we are going to verify in Part III below.

\medbreak By (6.4) for every $1\leq\nu\leq J$ there exists $m_\nu$ such that for a generic point $x^{(\nu-1)}$ of $X^{(\nu-1)}$ the curvature of the restriction, to $\left(\pi^{(\nu)}\right)^{-1}\left(x^{(\nu-1)}\right)$, of $R^0\left(\pi^{(\nu)}\right)_*\left(m_\nu K_{X^{(\nu)}|X^{(\nu-1)}}\right)$ is strictly positive at a generic point of $\left(\pi^{(\nu)}\right)^{-1}\left(x^{(\nu-1)}\right)$.  Here for the sake of simplicity of presentation we have suppressed some details by assuming that we have the simpler case of locally free $R^0\left(\pi^{(\nu)}\right)_*\left(m_\nu K_{X^{(\nu)}}\right)$ over $X^{(\nu-1)}$.  We have the following relations $K_{X^(\nu)}=K_{X^{(\nu)}|X^{(\nu-1)}}\left(\left(\pi^{(\nu)}\right)^*K_{X^{(\nu-1)}}\right)$ on $X^{(\nu)}$ for $0\leq\nu\leq J$.  Let $\hat m=\prod_{\nu=0}^J m_\nu$.  By the standard use of $L^2$ estimates of $\bar\partial$, we can now construct enough elements of $\Gamma\left(X, p\hat mK_X\right)$ to make $\dim_{\mathbb C}\Gamma\left(X, p\hat mK_X\right)$ increase to the order $p^{\kappa_{\rm num}(X)}$ as $p\to\infty$.  This finishes the confirmation of the abundance conjecture except the verification of the coincidence of the numerically trivial fibration for the canonical line bundle with the numerically trivial fibration for the canonical line bundle which is to be done in Part III below.

\bigbreak\noindent{\bf PART III. Coincidence of Numerically Trivial Foliation and Fibration for Canonical Bundle.}

\bigbreak\noindent{\bf\S7. Technique of Nevanlinna's First Main Theorem for Proof of Compactness of Leaves of Foliation.}
We are going to use the technique of Nevanlinna's First Main Theorem to show that the numerically trivial foliation and the numerically trivial fibration coincide for the canonical line bundle of a compact complex algebraic manifold $X$ of complex dimension $n$.  We will argue by contradiction.  By replacing $X$ by a generic fiber in the numerically trivial fibration for $K_X$ and after resolving singularities, we can assume that a generic leaf of the numerically trivial foliation of $X$ for $K_X$ is Zariski dense in $X$ and we will use this assumption to arrive at a contradiction.  For our arguments to arrive at a contradiction, in order to avoid the distraction of technical details which deal with singularities we will explain our arguments by imposing regularity conditions in some steps so that we can more focus on the essential points of our arguments.  First we introduce the notations we want from the theory of Nevanlinna by presenting very briefly the First Main Theorem of Nevanlinna in the form which is convenient for our application.

\bigbreak\noindent(7.1) {\it The First Main Theorem of Nevanlinna.}  The key technique to prove the coincidence of the numerically trivial foliation and fibration for the canonical line bundle is the method of integration by parts twice introduced by Nevanlinna for
the proof of the First Main Theorem of Nevanlinna theory [Nevanlinna1925, p.18, Erster Hauptsatz].

\medbreak For a meromorphic
section $s$ of a holomorphic line bundle $L$ with a smooth Hermitian metric $e^{-\psi}$ over a compact complex manifold $Y$, the element
of $H^1\left(Y,\Omega_Y^1\right)$ represented by the curvature form $\Theta_L:=\frac{\sqrt{-1}}{2\pi}\partial\bar\partial\psi$ of $e^{-\psi}$ agrees with the
element of $H^1\left(Y,\Omega_Y^1\right)$ represented by the divisor $Z_s$ of $s$ as a closed $(1,1)$-current.  This is simply the direct consequence of
$$\int_Y\frac{\sqrt{-1}}{2\pi}\partial\bar\partial\log\left\|s\right\|^2=0\leqno{(7.1.1)}$$
from Stokes's theorem, where $\left\|s\right\|^2=\left|s\right|^2e^{-\psi}$ is the pointwise square norm of $s$ with respect to the metric $e^{-\psi}$.

\medbreak When $Y$ is noncompact, the statement no longer holds in general, because in general an application of Stokes's theorem yields a ``boundary term'' (as the limit of the boundary for relatively compact subdomains of $Y$ exhausting $Y$)
which may not be zero.  The theory of Nevanlinna studies the special case of the noncompact manifold ${\mathbb C}$ when the line bundle and the section over ${\mathbb C}$ are obtained respectively by pulling back a line bundle $L$ on $Y$ with Hermitian metric $e^{-\psi}$ and a meromorphic section $s$ of $L$ on $Y$ via a holomorphic map $f:{\mathbb C}\to Y$.  In the theory of Nevanlinna the First Main Theorem is
$$
T\left(r,f,\Theta _L\right) = m\left(r,f,Zs\right) + N\left(r,f,Zs\right) + O(1)\ \ {\rm as\ \ }r\to\infty,\leqno{(7.1.2)}
$$
where
$$T\left(r,f,\Theta _L\right)=\int_{\rho=0}^r\frac{d\rho}{\rho}\int_{\Delta_\rho}f^*\Theta_L$$
is Nevanlinna's characteristic function of $f$,
$$m\left(r,f,Zs\right)=\frac{1}{2\pi}\int_{\theta=0}^{2\pi}\log\left(f^*\frac{1}{\left\|s\right\|^2}\right)$$ is the proximity function of the divisor $Zs$ for $f$, and
$$N\left(r,f,Zs\right)=\int_{\rho=0}^r\frac{n\left(\rho,f,Zs\right)d\rho}{\rho}$$ is the counting function of the divisor $Zs$ for $f$ with $n\left(\rho,f,Zs\right)$
being the number of zeroes of $f^*s$ minus the number of poles of $f^*s$ on $\Delta_\rho$ with multiplicities counted.    The method in the theory of Nevanlinna to derive (7.1.2) is to replace the one single application of Stokes's theorem in $(7.1.1)$ by using the following integration by parts twice
$$\int_{\rho=0}^r\frac{d\rho}{\rho}\int_{\Delta_\rho}\frac{\sqrt{-1}}{2\pi}\partial\bar\partial\log\left(f^*\left\|s\right\|^2\right)$$
to get $(7.1.2)$. When $Y$ is the Riemann sphere ${\mathbb P}_1$ and $L$ is the hyperplane section line bundle of ${\mathbb P}_1$ with the natural metric whose curvature $\Theta_L$ is the Fubini-Study metric of ${\mathbb P}_1$, the First Main Theorem is equivalent to the Poisson-Jensen formula (2.1.2.2) with $\log\left|f\left(re^{i\theta}\right)\right|$ rewritten as $$\log^+\left|f\left(re^{i\theta}\right)\right|-\log^+\frac{1}{\left|f\left(re^{i\theta}\right)\right|},$$ where $\log^+$ means the maximum of $\log$ and $0$.  We now express the First Main Theorem (7.1.2) in notations which are more suitable for our purpose.

\medbreak For a function or a (1,1)-form  $\eta$ on the topological closure of $\Delta_r$ we introduce the double integral  $${\mathcal I}_r(\eta ) =
\int ^r_{\rho =0} {d\rho \over \rho } \int _{\Delta_\rho}\eta.$$  For a function
$g$ on the boundary of $\Delta_r$ we introduce its average  $${\mathcal A}_r(g) = {1\over 2\pi } \int ^{2\pi }_{\theta =0}
g\left(r e^{i\theta }\right) d\theta.$$   The divergence theorem applied to
a function $g$ defined on the topological closure of $\Delta_\rho$ yields
$$
\int_{\Delta_\rho}\Delta g=\int ^{2\pi }_{\theta =0} \left({\partial \over \partial\rho} g\left(\rho e^{i\theta }\right)\right)\rho d\theta
=\rho\frac{d}{d\rho}\int ^{2\pi }_{\theta =0}g\left(\rho e^{i\theta }\right)d\theta\leqno{(7.1.3)}
$$
which yields $${\mathcal I}_r\left({1\over \pi }\frac{\partial^2 g}{\partial\zeta\partial\overline{\zeta}}
\right) = {1\over 2} {\mathcal A}_r(g) - {1\over 2} g(0)\leqno{(7.1.4)}$$
upon applying the operator $\int_{\rho=0}^r\frac{d\rho}{\rho}$ to both sides of $(7.1.3)$, where $\zeta$ is the coordinate of ${\mathbb c}$ and $\Delta  = 4\frac{\partial^2}{\partial\zeta\partial\overline{\zeta}}$ is the Laplace operator on ${\mathbb C}$.

\medbreak When the identity $(7.1.4)$ is applied to $g=f^*\left\|s\right\|^2$, from ${\sqrt{-1}\over 2\pi } \partial \bar \partial  \log  \|s\|^2
= - \Theta _L + Zs$ as closed $(1,1)$-currents on ${\mathbb C}$ we get
$$
{\mathcal I}_r\left(-f^*\left(\Theta _L + Zs\right)
\right) = {\mathcal A}_r\left(\log  \|f^*s\|\right) - \log  \|f^*s\|(0)
$$
and the First Main Theorem $(7.1.2)$
corresponds to the following rearrangement of terms
$$
{\mathcal I}_r\left(f^*\Theta _L\right) = {\mathcal A}_r\left(\log\frac{1}{\|f^*s\|}\right)+{\mathcal I}_r\left(f^*Zs\right)+\log\|f^*s\|(0)\leqno{(7.1.5)}
$$
term for term.

\bigbreak\noindent(7.2) {\it Additional Assumptions of Regularity.} For the sake of simplicity of presentation let us make the following three additional assumptions on the numerically trivial foliation ${\mathcal F}$ of the compact complex manifold $X$ of complex dimension $n$ for the canonical line bundle $K_X$.
\begin{itemize}\item[(i)] The leaves of ${\mathcal F}$ are regular of complex dimension $d$ and one of them $\hat M$ is dense in $X$.
\item[(ii)] The curvature current $\Xi_{K_X}$ of $K_X$ (which is defined in (1.1)) is smooth and its kernel defines precisely the smooth holomorphic foliation ${\mathcal F}$.
\item[(iii)]  The curvature form $\Xi_{K_X}$ is strictly positive on the normal bundle $N_{\mathcal F}$ of ${\mathcal F}$, which is defined as the quotient of $T_X$ by the tangent bundle $T_{\mathcal F}$, where $T_{\mathcal F}$ is the subbundle of $T_X$ consisting of all elements of $T_X$ tangential to the leaves of ${\mathcal F}$.
\end{itemize}
The general case with singularities necessitates tedious modifications, but the key arguments remain the same as for the case with additional assumptions of regularity.  First of all we introduce the following adjunction argument for foliations which is the natural extension of the usual simple adjunction formula for a complex submanifold relating the canonical line bundle of the submanifold, the restriction to the submanifold of the canonical line bundle of the ambient manifold and the normal bundle of the submanifold.

\bigbreak\noindent(7.3) {\it Adjunction of Foliation.}  Since $\Xi_{K_X}$ is the curvature form of $K_X$ and is assumed to be smooth, we can find a smooth volume form $\tau$ on $X$ such that 
$$\Xi_{K_X}=\frac{\sqrt{-1}}{2\pi}\partial\bar\partial\log\left|\tau\right|^2\quad{\rm on\ \ }X.\leqno{(7.3.1)}$$  We claim that the line bundle $\wedge^dT_{\mathcal F}$ is flat on $X$ along the directions of $T_{\mathcal F}$ in the sense that there exists a metric of $\wedge^dT_{\mathcal F}$ whose curvature when restricted to $T_{\mathcal F}$ is identically zero.  The verification of the claim is as follows.

\medbreak Let $e^{-\varphi}=\frac{1}{\left|\tau\right|^2}$.  The normal bundle $N_{\mathcal F}$ of the foliation ${\mathcal F}$ (which is a vector bundle of rank $n-d$ over $X$) is given by $N_{\mathcal F}=T_X\left/T_{\mathcal F}\right.$. Since $\Xi_{K_X}$ is semipositive, for every $\xi\in T_{\mathcal F}$ and $\eta\in T_X$ not only we have $\Xi_{K_X}\left(\xi,\bar\xi\right)=0$ from $\xi\in T_{\mathcal F}$ by definition of ${\mathcal F}$ and $T_{\mathcal F}$, but we have also $\Xi_{K_X}\left(\xi,\bar\eta\right)=0$ from the semipositivity of $\Xi_{K_X}$.  As a result, at any point $P$ of $X$ with a sufficiently small coordinate polydisk neighborhood $U$ of $P$, not only on $U$ the restriction of $\varphi$ to a leaf of the foliation ${\mathcal F}$ is pluriharmonic, but there exists a holomorphic function $H$ on $U$ such that on $U$ the restriction of $\varphi$ to a leaf of the foliation ${\mathcal F}$ is equal to the restriction of the real part ${\rm Re}\,H$ of $H$ to that leaf of the foliation ${\mathcal F}$.  With the use of a new local trivialization defined by $e^{\frac{H}{2}}$ on $U$ with respect to the old local trivialization, we can assume that $\varphi$ is constant along the leaf of ${\mathcal F}$ on $U$.

\medbreak We have a holomorphic map $\pi$ from $U$ to an open polydisk $W$ in ${\mathbb C}^{n-d}$ such that its fibers $\pi^{-1}(w)$ for $w\in W$ is an open subset of a leaf of the foliation ${\mathcal F}$.  Since $\varphi$ is locally constant on the leaves of ${\mathcal F}$ on $U$ there exists a smooth function $\psi$ on $W$ such that $\varphi|_U=\psi\circ\pi$.   Let $\theta=\frac{\sqrt{-1}}{2\pi}\partial\bar\partial\psi$ on $W$. Then $\Xi_{K_X}=\pi^*\theta$ on $U$.

\medbreak We consider $\left(\Xi_{K_X}|_{N_{\mathcal F}}\right)^{n-d}$ which on $U$ is equal to $\left(\pi^*\theta\right)^{n-d}$.  Upon taking $\partial\bar\partial\log$ we conclude that
$$\partial\bar\partial\log\left(\left(\Xi_{K_X}|_{N_{\mathcal F}}\right)^{n-d}\right)=\pi^*\left(\partial\bar\partial\log\left(\theta^{n-d}\right)\right)\quad{\rm on\ }U.\leqno{(7.3.2)}$$
Let $v$ be a smooth positive $(d,d)$-form along the leaves of ${\mathcal F}$ so that $v$ defines a metric for the holomorphic line bundle $\wedge^d T_{\mathcal F}$ of $X$.  Then $v\wedge\left(\Xi_{K_X}\right)^{n-d}$ is a well-defined global positive definite $(n,n)$-form on $X$ and we can write $v\wedge\left(\Xi_{K_X}\right)^{n-d}=F\tau$ for some positive-valued function $F$ on $X$.  By replacing $v$ by $\frac{v}{F}$, we can assume that we have chosen $v$ such that $v\wedge\left(\Xi_{K_X}\right)^{n-d}=\tau$.  Taking $\frac{\sqrt{-1}}{2\pi}\partial\bar\partial\log$ of both sides of $v\wedge\left(\Xi_{K_X}\right)^{n-d}=\tau$, we conclude from (7.3.1) and (7.3.2) that the curvature $-\frac{\sqrt{-1}}{2\pi}\partial\bar\partial\log v$ of the metric $v$ of $\wedge^d T_{\mathcal F}$ vanishes when restricted to $T_{\mathcal F}$.  This finishes the verification of the claim.

\bigbreak\noindent(7.4) {\it Special Case of Leaf Dimension One.}  We first consider the simpler case where the leaf dimension $d$ is $1$.  For the case of a general dimension $d$ the argument is the same, but it is easier to first explain the argument in the case of $d=1$ without the distraction of technical complications of the case of a general dimension $d$.  We now start our argument aiming to derive a contradiction by assuming the leaf dimension $d$ to be $1$ and assuming the existence of a dense leaf $\hat M$.   The universal cover $M$ of the leaf $\hat M$ is simply biholomorphic to ${\mathbb C}$ for the following reason.   Since $K_X$ is flat along the leaves of ${\mathcal F}$, we can find a metric of $T_{\mathcal F}$ with zero curvature.  This means that the K\"ahler metric on the simply connected Riemann surface $M$ induced from the flat metric of $T_{\mathcal F}$ is complete with zero curvature.  Hence $M$ is biholomorphic to ${\mathbb C}$.  Denote by $f:{\mathbb C}\to\hat M$ the universal covering map for the noncompact complex Riemann surface $\hat M$.

\bigbreak\noindent(7.5) {\it Application of First Main Theorem to 1-Dimensional Leaf of Numerically Trivial Foliation.}  We now have the universal cover $f:{\mathbb C}\to\hat M$ of the dense leaf $\hat M$ of complex dimension $1$ in the numerically trivial foliation ${\mathcal F}$ for the canonical line bundle $K_X$ of the compact complex algebraic manifold $X$ of complex dimension $n$.  Fix a point $\hat Q$ in $\mathbb C$ with $Q=f\left(\hat Q\right)\in\hat M$.  To explain the argument, we assume that $\dim_{\mathbb C}\Gamma\left(X, mK_X+A\right)\geq\left\lfloor\alpha m\right\rfloor$ for some $\alpha>0$ and for all $m\in{\mathbb N}$.  Note that when $\alpha>0$ is small, the statement $\dim_{\mathbb C}\Gamma\left(X, mK_X+A\right)\geq\left\lfloor\alpha m\right\rfloor$ for $1\leq m<\frac{1}{\alpha}$ simply means the trivial statement  $\dim_{\mathbb C}\Gamma\left(X, mK_X+A\right)\geq 0$.  Choose a nonzero element $s_m$ of $\Gamma\left(X, mK_X+A\right)$ such that the restriction of $s_m$ to $\hat M$ vanishes to order $\geq\frac{1}{2}\left\lfloor\alpha m\right\rfloor$ at $Q$ and we normalize $s_m$ by multiplying it by a positive constant so that $\max_{P\in X}\left\|s_m\left(P\right)\right\|=1$ and $\left\|s_m\left(P_m\right)\right\|=1$ at some point $P_m\in X$ (which may not be on $\hat M$).  Here $\left\|s_m\right\|^2$ means $\left|s_m\right|^2\,e^{-m\varphi_{K_X}-\varphi_A}$, where $e^{-\varphi_A}$ is a smooth metric for $A$ whose curvature form is positive definite and $e^{-\varphi_{K_X}}$ is the metric of $K_X$ whose curvature current is $\Xi_{K_X}$ from (1.1) which is assumed smooth and satisfies (ii) and (iii) of (7.2).  Note that since $\hat M$ is dense in $X$, we know that the restriction of $s_m$ to $\hat M$ is not identically zero on $\hat M$.

\medbreak We choose some subsequence $m_\nu$ such that $P_{m_\nu}\to P_*$ in $X$ as $\nu\to\infty$ for some $P_*\in X$.  Choose a local coordinate system $\left(z_1,\cdots,z_n\right)$ of $X$ centered at $P_*$ and denote by $B_{r,P}$ the coordinate ball centered at $P$ of radius $r$ with respect to $\left(z_1,\cdots,z_n\right)$.  Choose trivializations of $K_X$ and $A$ on an open neighborhood $U_{P_*}$ of $P_*$ in $X$ so that any holomorphic section of $mK_X+A$ on any open subset $W$ of $U_{P_*}$ can be represented by a holomorphic function on $W$ according to the chosen trivializations of $K_X$ and $A$ on $U_{P_*}$.   Moreover, we choose the trivialization of $K_X$ on $U_{P_*}$ with the additional property that with respect to the chosen trivialization of $K_X$ on $U_{P_*}$ the function $\varphi_{K_X}$ as well as all its first-order derivatives vanishes at $P_*$.

\medbreak Let $r_{m_\nu}>0$ be the supremum of $r>0$ such that $s_{m_\nu}$ is nowhere zero on $B_{r,P_{m_\nu}}$.  Denote by $\left(s_{m_\nu}\right)^{\frac{1}{m_\nu}}$ a branch of the holomorphic function $s_{m_\nu}$ on $B_{r_{m_\nu},P_{m_\nu}}$.  Since
$\left|\left(s_{m_\nu}\right)^{\frac{1}{m_\nu}}\right|^2e^{-\varphi_{K_X}-\frac{1}{m_\nu}\varphi_A}=\left\|\left(s_{m_\nu}\right)^{\frac{1}{m_\nu}}\right\|^2\leq 1$ on $B_{r_\nu,P_{m_\nu}}$ and since both $\varphi_{K_X}$ and $\varphi_A$ are smooth on $X$, by applying Cauchy's integral formula for the first-order derivatives of the branch of $\left(s_{m_\nu}\right)^{\frac{1}{m_\nu}}$ and using $\left\|s_{m_\nu}\left(P_{m_\nu}\right)\right\|=1$, we conclude that there exist $\nu_0\in{\mathbb N}$ and $r_*>0$ such that $r_{m_\nu}\geq r_*$ for $\nu\geq\nu_0$ and
$\left|\left(s_{m_\nu}\right)^{\frac{1}{m_\nu}}\right|\geq\frac{1}{2}$ on $B_{r_*,P_{m_\nu}}$ for $\nu\geq\nu_0$.   Moreover, since the function $\varphi_{K_X}$ as well as all its first-order derivatives vanishes at $P_*$, we can assume that $\nu_0\in{\mathbb N}$ is chosen so large and $r_*>0$ is chosen so small that $\left\|\left(s_{m_\nu}\right)^{\frac{1}{m_\nu}}\right\|\geq\frac{1}{3}$ on $B_{r_*,P_{m_\nu}}$ for $\nu\geq\nu_0$.

\medbreak Since $\hat M$ is dense in $X$ and since
$r_{m_\nu}\geq r_*$ for $\nu\geq\nu_0$ and $P_{m_\nu}\to P_*$ as $\nu\to\infty$,
it follows
that there exist some integer $\nu_1\geq\nu_0$ and some point $P_{\#}\in\hat M$ such that $P_{\#}\in B_{\frac{r_*}{2},P_{m_\nu}}$ for $\nu\geq\nu_1$.  Thus  $\left\|\left(s_{m_\nu}\right)^{\frac{1}{m_\nu}}\left(P_{\#}\right)\right\|\geq\frac{1}{3}$ for $\nu\geq\nu_1$.   Let $\hat P\in{\mathbb C}$ such that $P_{\#}=f\left(\hat P\right)$.  By using a translation of the coordinate of ${\mathbb C}$, we can assume without loss of generality that $\hat P$ is the origin of ${\mathbb C}$.  Now we apply $(7.1.5)$ to $L=m_\nu K_X+A$ and $s=s_{m_\nu}$ for $\nu\geq\nu_1$ and divide the result by $m_\nu$ to get
$$
\displaylines{(7.5.1)_\nu\qquad\qquad
{\mathcal I}_r\left(\frac{1}{m_\nu}f^*\left(m_\nu\Theta_{K_X}+\Theta_A\right)\right)\hfill\cr\hfill =
\frac{1}{m_\nu}{\mathcal A}_r\left(\log\frac{1}{\|f^*s_{m_\nu}\|}\right)+{\mathcal I}_r\left(\frac{1}{m_\nu}f^*Zs_{m_\nu}\right)+\frac{1}{m_\nu}\log\|f^*s_{m_\nu}\|(0).\cr}
$$
Let $\hat r$ be the distance of $\hat Q$ from the origin in ${\mathbb C}$.  Take an arbitrary $r\geq 2\hat r$. Since $\left\|s_{m_\nu}\right\|\leq 1$ on $X$, it follows that ${\mathcal A}_r\left(\log\frac{1}{\|f^*s_{m_\nu}\|}\right)$ is nonnegative.  From $\left\|\left(s_{m_\nu}\right)^{\frac{1}{m_\nu}}\left(P_{\#}\right)\right\|\geq\frac{1}{3}$ for $\nu\geq\nu_1$ and $P_{\#}=f(0)$ it follows that $\frac{-1}{m_\nu}\log\|f^*s_{m_\nu}\|(0)\leq\log 3$ for $\nu\geq\nu_1$.  Since the restriction of $s_m$ to $\hat M$ vanishes to order $\geq\frac{1}{2}\left\lfloor\alpha m_\nu\right\rfloor$ at $Q$, it follows that there exists some $C_*$ independent of $r\geq 2\hat r$ such that
$${\mathcal I}_r\left(\frac{1}{m_\nu}f^*Zs_{m_\nu}\right)\geq \frac{1}{3}\left\lfloor\alpha\right\rfloor\,\log r-C_*$$
for all $\nu\geq\nu_1$.  Since the pullback of $\Theta_{K_X}$ to $\hat M$ is identically zero, we have
$$
\lim_{\nu\to\infty}{\mathcal I}_r\left(\frac{1}{m_\nu}f^*\left(m_\nu\Theta_{K_X}+\Theta_A\right)\right)=0
$$
for any fixed $r\geq 2\hat r$.  From $(7.5.1)_\nu$ as $\nu\to\infty$ we obtain the inequality $\frac{1}{3}\left\lfloor\alpha\right\rfloor\,\log r-C_*\leq\log 3$, which would give a contradiction for a sufficiently large $r\geq 2\hat r$.  For the case at hand this finishes the proof of the compactness of the leaves of the numerically trivial foliation for the canonical line bundle.

\bigbreak\noindent(7.6) {\it Comparison with Proof of Bombieri's Higher-Dimensional Formulation of Technique of Gelfond-Schneider.}  The argument in (7.1), (7.4), and (7.5) for the case of leaf dimension $1$ follows the same line as the proof of Bombieri's higher-dimensional formulation of the technique of Gelfond-Schneider given in (2.1).  In (7.5) only an inequality from the First Main Theorem is used, which is analogous to the inequality (2.1.2.3) in the Gelfond-Schneider-Bombieri argument.  However, there are two differences.

\medbreak\noindent(7.6.1) The first difference is the way to handle the lower bound of the $O(1)$ term.  In the Gelfond-Schneider-Bombieri argument the lower bound of $\log\left|f\left(w_0\right)\right|$ is obtained by differentiating (2.1.1.1) $s-1$ times and from the use of the chain rule and the height estimates of
${\rm size}\left(a_{j_1,\cdots,j_{d+1}}\right)=O\left(L\right)$.  On the other hand, in the the argument of the case of leaf dimension $1$ the lower bound comes from the normalization $\max_{P\in X}\left\|s_m\left(P\right)\right\|=1$ and $\left\|s_m\left(P_m\right)\right\|=1$ and the use of Cauchy's integral formula for the first-order derivatives.

\medbreak\noindent(7.6.2) The second difference is the way in which the function $P\left(f_1,\cdots,f_{d+1}\right)$ and the section $s_m$ with vanishing order conditions are produced.  The construction of $P\left(f_1,\cdots,f_{d+1}\right)$ in the Gelfond-Schneider-Bombieri technique is more complicated and requires the algebraic independence of $f_1,\cdots,f_{d+1}$ over the number field $K$, because the height estimate of ${\rm size}\left(a_{j_1,\cdots,j_{d+1}}\right)=O\left(L\right)$ for its coefficients requires the degree of $P\left(\xi_1,\cdots,\xi_{d+1}\right)$ in each the $d+1$ variables to be $\leq J$ which is lower than the lower bound $L$ of the vanishing order of $P\left(f_1,\cdots,f_{d+1}\right)$ at points of $S_m$ in a specified way.  On the other hand the construction of the section $s_m$ in the argument of the case of leaf dimension $1$ is a simple argument involving the growth of $\dim_{\mathbb C}\Gamma\left(X,mK_X+A\right)$ to be at least of linear order in $m$.

\bigbreak\noindent(7.7) {\it General Case of Higher Leaf Dimension.}  We will only give a brief sketch of the modifications of (7.1), (7.4), and (7.5) needed to handle the general case of higher leaf dimension.  Again we start out with a dense leaf of complex dimension $d\geq 2$ and argue to arrive at a contradiction.  We use the additional assumptions of regularity listed in (7.2) to suppress the distracting details arising from singularities.  First let us remark that the First Main Theorem for holomorphic maps from ${\mathbb C}$ to a complex manifold $Y$ can be interpreted as Green's formula for Green's function on a disk in ${\mathbb C}$ in the following way.

\bigbreak\noindent(7.7.1) {\it First Main Theorem as Green's Formula of Green's Function.}  Let $M$ be a complete K\"ahler manifold with K\"ahler form $\omega_M$ and let $P_0$ be a point of $M$.  Let $r$ be the distance function on $M$ measured from $P_0$.  For $R>0$ let $B_R$ be the ball of radius $R$ centered at $P_0$ consisting of all points $P$ of $M$ whose distance $r(P)$ from $P_0$ with respect to $\omega_M$ is $<R$.  Let $G_R(P)$ (as a function of $P$) be Green's function for the domain $B_R$ and the point $P_0$ in the sense that the Laplacian of $G_R(P)$ (as a function of $P$) evaluated at $P_0$ is equal to Dirac delta at $P_0$ and the boundary value of $G_R(P)$ at $\partial B_R$ is identically zero.  For any function $F$ on the topological closure of $B_R$ Green's formula
$$
\int_{B_R}\left(F\Delta G_R-G_R\Delta F\right)=\int_{\partial B_R}\left(F\frac{\partial G_R}{\partial r}-G_R\frac{\partial F}{\partial r}\right),
$$
yields the following Green's function formula
$$
F\left(P_0\right)=\int_{\partial B_R} F\frac{\partial G_R}{\partial r}+\int_{B_R} G_R\Delta F.\leqno{(7.7.1.1)}
$$
The equation (2.1.2.2) is simply the special case of the Green's function formula $(7.7.1.1)$ with $M={\mathbb C}$ and $B_R=\Delta_R$ with $R=r$ and $F$ being the logarithm of the absolute value of a meromorphic function on the topological closure of $\overline{B_R}$, because Green's function $G_R(\zeta)$ for $B_R$ is given by $\frac{1}{4\pi}\log\frac{|\zeta|}{R}$ with $\zeta$ being the coordinate of ${\mathbb C}$.  When we have a holomorphic map $f$ from ${\mathbb C}$ to a complex complex manifold $Y$ and a holomorphic line bundle $L$ with Hermitian metric $e^{-\psi}$ over $Y$ and a holomorphic section $s$ of $L$ over $Y$, the First Main Theorem is simply the special case of $(7.7.1.1)$ with $M={\mathbb C}$ and $B_R=\Delta_R$ with $R=r$ and $F=\log\left(\left|s\right|^2 e^{\psi}\right)$.

\bigbreak\noindent(7.7.2) {\it Replacement in the Case of General Leaf Dimension for ${\mathbb C}$ as the Universal Cover of 1-Dimensional Leaf.}  When the leaf dimension $d$ is $1$, we have the universal covering map $f:{\mathbb C}\to\hat M$ for a dense leaf $\hat M$ of complex dimension $1$ in $X$.   Now we are in the case of general leaf dimension $d\geq 2$, we have to replace the use of the universal covering map $f:{\mathbb C}\to\hat M$ of the leaf $\hat M$ by a more complicated construction.  If we have a universal covering map ${\mathbb C}^d\to\hat M$ for our $d$-dimensional leaf $\hat M$, the same argument which we use for the $1$-dimensional leaf would go through in the same way.  In general we do not have this ideal situation, but we want to get close enough to the ideal situation for the argument to work.  That is the reason for the more complicated construction given immediately below.

\bigbreak\noindent(7.7.2.1) {\it Ricci Flat Metric Along Leaves of Foliation.} First of all, from the flatness of the canonical line bundle $\wedge^d\left(T_{\mathcal F}\right)^*$ of the leaves of ${\mathcal F}$ we construct a K\"ahler metric for $T_{\mathcal F}$ with zero Ricci curvature in the sense that there is a positive-definite $(1,1)$-form $\theta$ along the leaves of ${\mathcal F}$ which defines a smooth metric for the vector bundle $T_{\mathcal F}$ of rank $d$ on $X$ such that $\partial\bar\partial\log\left(\wedge^d\theta\right)$ is identically zero when restricted to $T_{\mathcal F}$.

\medbreak The construction of such a K\"ahler metric with zero Ricci curvature along the leaves of the foliation is obtained by modifying Yau's construction of K\"ahler-Einstein metrics for a compact complex K\"ahler manifold $Y$ of complex dimension $\ell$ with trivial canonical line bundle [Yau1978].  One key modification is that the integration of a smooth $(\ell,\ell)$-form on $Y$ in Yau's construction is replaced by the integration of the exterior product of a smooth $(d,d)$-form with $\left(\Xi_{K_X}\right)^{n-d}$ over $X$ in our case.

\bigbreak\noindent(7.7.2.2) {\it Partial Removal of Cut-Locus.}  Let $\sigma:{\mathbb C}^d\to\hat M$ be the (smooth and non-holomorphic) exponential map from a point $P_0$ of $\hat M$ with respect to the complete K\"ahler metric of $\hat M$ of zero Ricci curvature which is induced from the K\"ahler metric for $T_{\mathcal F}$ with zero Ricci curvature.  We give ${\mathbb C}^d$ (outside the zero-set $Z$ of the differential $d\sigma$ of $\sigma$) the complex structure and the K\"ahler metric from $\hat M$ induced by the smooth map $\sigma$ and we use $M$ to denote ${\mathbb C}^d$ with this new complex structure and the new K\"ahler metric (outside $Z$).  This map is introduced to partially handle the cut-locus of $\hat M$ for its K\"ahler-Einstein metric from the point $P_0$, namely the part which is not due to the vanishing of the differential $d\sigma$ of the exponential map $\sigma$ but due to the intersection of two different smooth geodesics both emanating from $P_0$.  We call $P_0$ the origin of $M$.

\bigbreak\noindent(7.7.2.3) There are two reasons for introducing $M$.  They concern the Green's function formula (7.7.1.1) of $M$.  The first reason is the special properties of Green's function on a complete K\"ahler manifold of zero Ricci curvature obtained from the Jacobi field equation along a geodesic ({\it cf.,} [Berger-Gauduchon-Mazet1971, Cheeger-Ebin1975, Siu-Yau1977]).  These special properties of Green's function make the formula (7.7.1.1) close enough to the formula in the special case of the complex manifold being ${\mathbb C}^n$ so that our argument for the use of the formula (7.7.1.1) goes through like the case of ${\mathbb C}^n$.  Though there is an exceptional set $Z$ in the manifold $M$ introduced here, the argument which depends on these special properties of Green's function still works.  

\medbreak The second reason is that, with respect to a K\"ahler metric of a complex manifold of complex dimension $n$ with K\"ahler form $\omega$, the Laplacian $\Delta F$ is given by
$$
\left(\sqrt{-1}\,\partial\bar\partial F\right)\wedge\frac{\omega^{n-1}}{(n-1)!}
=\left(\Delta F\right)\frac{\omega^n}{n!}
$$
so that, when it is applied to the case of $F=\left|s\right|^2e^{-\psi}$ with $s$ holomorphic and $e^{-\psi}$ being the metric of a line bundle, it naturally becomes the exterior product of $\omega^{n-1}$ with the difference of two closed positive $(1,1)$-currents and we can get the estimates we want, as in (7.1) and (7.4).

\bigbreak\noindent(7.7.3) {\it Normalized Amply Twisted Pluricanonical Section with Appropriate Vanishing Order at Prescribed Point.}  As in the argument for the case of leaf dimension $1$, we choose $Q\in\hat M$ and then choose $s_m\in\Gamma\left(X,mK_X+A\right)$ such that the restriction $s_m|_{\hat M}$ of $s_m$ to $\hat M$ vanishes at $Q$ to an order $\geq\left\lfloor\alpha m\right\rfloor$ at $Q$ (where $\alpha$ is a positive number as in (7.5)).  We then normalize $s_m$ so that $\left\|s_m\right\|^2=\left|s_m\right|^2e^{-m\varphi_{K_X}-\varphi_A}\leq 1$ everywhere on $X$ with $\left\|s_m\right\|^2=\left|s_m\right|^2e^{-m\varphi_{K_X}-\varphi_A}=1$ at some point $P_m$ of $X$.  Let $Zs_m$ be the divisor of $s_m$ in $X$ and $\tilde Z_m$ be its pullback in $M$ via the composite map $M\to\hat M\hookrightarrow X$.  As in (7.5) we can assume that the origin in $M$ is mapped to a point $P_{\#}$ of $\hat M$ such that $\left\|s_{m_\nu}\right\|\left(P_{\#}\right)\geq c_1$ for some positive number $c_1$ and some subsequence $m_\nu$ of $m$.  By applying the Green's function formula (7.7.1.1) to $M$ as explained in (7.7.2.3) and using the balls $B_r$ in $M$ centered at the origin and of radius $r$, we obtain as in (7.5) positive numbers $c_2, C_3, R_*$ such that
$$
c_2\log R\leq \lim_{\nu\to\infty}\frac{1}{m_\nu}\int_{r=0}^R\left({\rm volume\ of\ }B_r\cap\tilde Z_{m_\nu}\right)\,\frac{dr}{r} \leq C_3
$$
for $R\geq R_*$ sufficiently large.  This gives the contradiction as $R\to\infty$.

\bigbreak\noindent{\it References}

\medbreak\noindent[Bauer-et-al2002] Thomas Bauer, Fr\'ed\'eric Campana, Thomas Eckl, Stefan Kebekus, Thomas Peternell, S\l awomir Rams, Tomasz Szemberg, and Lorenz Wotzlaw,
A reduction map for nef line bundles. {\it Complex geometry} (G\"ottingen, 2000), 27--36, Springer, Berlin, 2002.

\medbreak\noindent[Berger-Gauduchon-Mazet1971]
Marcel Berger, Paul Gauduchon and Edmond Mazet,
{\it Le spectre d'une vari\'et\'e riemannienne}.
Lecture Notes in Mathematics, Vol. \textbf{194} Springer-Verlag, Berlin-New York 1971.

\medbreak\noindent
[Berndtsson2006] Bo Berndtsson, Positivity of direct image bundles and convexity on the space of K\"ahler metrics
math arXiv: math.CV/0608385.

\medbreak\noindent
[Berndtsson2009] Bo Berndtsson,
Curvature of vector bundles associated to holomorphic fibrations.
{\it Ann. of Math.} \textbf{169} (2009), 531--560.

\medbreak\noindent[Berndtsson-Paun2008]
Bo Berndtsson and Mihai Paun, Bergman kernels and the pseudoeffectivity of relative canonical bundles. arXiv:math/0703344.
{\it Duke Math. J.} \textbf{145} (2008), 341--378.

\medbreak\noindent[Bombieri1970] Enrico Bombieri, Algebraic values of meromorphic maps, {\it Invent.
Math.} \textbf{10} (1970), 267-287.  Addendum. {\it Invent.
Math.} \textbf{11} (1970), 163--166.

\medbreak\noindent[Bombieri-Lang1970] Enrico Bombieri and Serge Lang,
Analytic subgroups of group varieties.
{\it Invent. Math.} \textbf{11} (1970), 1--14.

\medbreak\noindent[Brieskorn1970]
Egbert Brieskorn,
Die Monodromie der isolierten Singularit\"aten von Hyperfl\"achen.
{\it Manuscripta Math.} \textbf{2} (1970), 103--161.

\medbreak\noindent
[Budur2009] Nero Budur, Unitary local systems, multiplier ideals, and polynomial
periodicity of Hodge numbers. {\it Adv. Math.} \textbf{221} (2009), 217--250.

\medbreak\noindent[Campana-Peternell-Toma2007]
Frederic Campana, Thomas Peternell, and Matei Toma
Geometric stability of the cotangent bundle and the universal cover of a projective manifold. arXiv:math/0405093.

\medbreak\noindent[Cheeger-Ebin1975] Jeff Cheeger and David G. Ebin,
{\it Comparison theorems in Riemannian geometry}.
North-Holland Mathematical Library, Vol. \textbf{9}. North-Holland Publishing Co., Amsterdam-Oxford; American Elsevier Publishing Co., Inc., New York, 1975.

\medbreak\noindent[Deligne-Mostow1986] Pierre Deligne and George Daniel Mostow,
Monodromy of hypergeometric functions and nonlattice integral monodromy.
{\it Inst. Hautes \'Etudes Sci. Publ. Math.} \textbf{63} (1986), 5--89.

\medbreak\noindent[Eckl2004a] Thomas Eckl, Tsuji's numerical trivial fibrations.
{\it J. Algebraic Geom.} \textbf{13} (2004), 617--639.

\medbreak\noindent[Eckl2004b] Thomas Eckl,
Numerically trivial foliations.
{\it Ann. Inst. Fourier (Grenoble)} \textbf{54} (2004), 887--938.

\medbreak\noindent[Euler1778]  Leonhard Euler, Specimen transformationi singularis serierum, Sept. 3, 1778, {\it Nova Acta Petropolitana} \textbf{7} (1801), 58-78.

\medbreak\noindent[Gauss1812] Johann Carl Friedrich Gauss, Disquisitiones generales circa seriem infinitam, Werke Band III, p.127, Formel I -- V.

\medbreak\noindent[Gelfond1934] A. O. Gelfond, Sur le septi\`eme Probl\`eme de D. Hilbert. {\it Comptes Rendus Acad. Sci. URSS Moscou} \textbf{2} (1934), 1-6. {\it Bull. Acad. Sci. URSS Leningrade} \textbf{7} (1934), 623-634.

\medbreak\noindent[Griffiths1968]
Phillip A. Griffiths,
Periods of integrals on algebraic manifolds. I. Construction and properties of the modular varieties.
{\it Amer. J. Math.} \textbf{90} (1968), 568--626.
II. Local study of the period mapping.
{\it Amer. J. Math.} \textbf{90} (1968), 805--865.

\medbreak\noindent[Hilbert1900] David Hilbert, Mathematische Probleme.  {\it Nachr. K\"onigl. Ges. der Wiss. zu G\"ottingen, Math.-Phys. Klasse} (1900), 251--297.

\medbreak\noindent[Iitaka1971] S. Iitaka, On $D$-dimensions of algebraic varieties.  {\it J. Math. Soc. Japan} \textbf{23} (1971), 356--373.

\medbreak\noindent
[Kawamata1982] Yujiro Kawamata, Kodaira dimension of algebraic fiber spaces over curves. {\it Invent.Math.} \textbf{66} (1982),
57–71.

\medbreak\noindent[Kawamata1985] Yujiro Kawamata, Pluricanonical systems on minimal algebraic varieties. {\it Invent. Math.} \textbf{79} (1985), 567--588.

\medbreak\noindent[Kodaira-Spencer1960] K. Kodaira and D. C. Spencer,
On deformations of complex analytic structures, III.
Stability theorems for complex structures. {\it Ann. of Math.} \textbf{71} (1960), 43-76.

\medbreak\noindent[Koll\'ar-Mori1998] J\'anos Koll\'ar and Shigefumi Mori,
{\it Birational geometry of algebraic varieties.}
With the collaboration of C. H. Clemens and A. Corti. Cambridge Tracts in Mathematics \textbf{134}. Cambridge University Press, Cambridge, 1998.

\medbreak\noindent[Lang1962] Serge Lang,
Transcendental points on group varieties.
{\it Topology} \textbf{1} (1962), 313--318.

\medbreak\noindent[Lang1965] Serge Lang,
Algebraic values of meromorphic functions.
{\it Topology} \textbf{3} (1965), 183--191.

\medbreak\noindent[Lang1966] Serge Lang,
{\it Introduction to transcendental numbers}. Addison-Wesley Publishing Co., Reading, Mass.-London-Don Mills, Ont. 1966.

\medbreak\noindent[LeVavasseur1893] Raymond Le Vavasseur,
Sur le syst\`eme d'\'equations aux d\'eriv\'ees partielles simultan\'ees auxquelles satisfait la s\'erie hyperg\'eométrique à deux variables, {\it Annales de la Facult\'e des Sciences de Toulouse} \textbf{7} (1893), 1-120.

\medbreak\noindent[Mostow1987] George Daniel Mostow, Braids, hypergeometric functions, and lattices, Bull. Amer. Math. Soc. 16 (1987), 225-246.

\medbreak\noindent[Nevanlinna1925]
Rolf Nevanlinna,
Zur Theorie der Meromorphen Funktionen.
{\it Acta Math.}
\textbf{46} (1925), 1--99.

\medbreak\noindent[Ohsawa-Takegoshi1987]
Takeo Ohsawa and Kensho Takegoshi,
On the extension of $L^2$ holomorphic functions.
{\it Math. Zeitschr.} \textbf{195} (1987),  197--204.

\medbreak\noindent[Paun2007]
Mihai Paun,
Siu's invariance of plurigenera: a one-tower proof.
{\it J. Differential Geom.} \textbf{76} (2007), 485--493.

\medbreak\noindent[Picard1885] \'Emile Picard, Sur les fonctions hyperfuchsiennes provenant des s\'eries hyperg\'eometrique de deux variables, {\it Ann. \"Ecole Norm. Sup.} \textbf{2} (1885), 357-384.

\medbreak\noindent[Nakayama2004] Noboru Nakayama,
{\it Zariski-decomposition and abundance.}
Mathematical Society of Japan Memoirs \textbf{14}, 2004.

\medbreak\noindent[Schneider1934] T. Schneider, Transzendenzuntersuchungen periodischer Funktionen. I, II. {\it J. reine angew. Math.} \textbf{172} (1934), 65-74.

\medbreak\noindent[Schumacher2009] Georg Schumacher, Positivity of relative canonical bundles for families of canonically polarized manifolds (2009), arXiv:0808.3259v2.

\medbreak\noindent[Schwarz1873]  Karl Hermann Amandus Schwarz, \"Uber diejenige F\"alle in welchen die Gaussische hypergeometrische Reihe eine algebraische Function ihres viertes elementes darstellt, {\it Crelle's J.} \textbf{75} (1873), 292-335.

\medbreak\noindent[Siegel1958] Carl L. Siegel,
On meromorphic functions of several variables.
{\it Bull. Calcutta Math. Soc.} \textbf{50} (1958), 165--168.

\medbreak\noindent[Simpson1993] Carlos Simpson,
Subspaces of moduli spaces of rank one local systems.
{\it Ann. Sci. \'Ecole Norm. Sup.} \textbf{26} (1993), 361--401.

\medbreak\noindent[Siu1974] Yum-Tong Siu,
Analyticity of sets associated to Lelong numbers and the extension of closed positive currents.
{\it Invent. Math.} \textbf{27} (1974), 53--156.

\medbreak\noindent[Siu1998] Yum-Tong Siu, Invariance of plurigenera.
{\it Invent. Math.} \textbf{134} (1998), 661--673.

\medbreak\noindent[Siu2002] Yum-Tong Siu,
Extension of twisted pluricanonical sections with plurisubharmonic weight and invariance of semipositively twisted plurigenera for manifolds not necessarily of general type. {\it Complex geometry} (G\"ottingen, 2000), 223--277, Springer, Berlin, 2002.

\medbreak\noindent[Siu2004] Yum-Tong Siu,
Invariance of plurigenera and torsion-freeness of direct image sheaves of pluricanonical bundles. {\it Finite or infinite dimensional complex analysis and applications}, 45--83,
{\it Adv. Complex Anal. Appl.}, \textbf{2}, Kluwer Acad. Publ., Dordrecht, 2004.

\medbreak\noindent[Siu2006] Yum-Tong Siu, A general non-vanishing
theorem and an analytic proof of the finite generation of the
canonical ring, arXiv:math/0610740.

\medbreak\noindent[Siu2008] Yum-Tong Siu,
Finite Generation of Canonical Ring by Analytic Method, arXiv:0803.2454,
{\it J. Sci. China} \textbf{51} (2008), 481-502.

\medbreak\noindent[Siu2009] Yum-Tong Siu,
Techniques for the Analytic Proof of the Finite Generation of the Canonical Ring, arXiv:0811.1211,
{\it Current Developments in Mathematics 2007}, ed. David Jerison et al, International Press 2009, pp.177-220.

\medbreak\noindent[Siu-Yau1977]
Yum-Tong Siu and Shing-Tung Yau,
Complete K\"ahler manifolds with nonpositive curvature of faster than quadratic decay.
{\it Ann. of Math.} \textbf{105} (1977), 225--264. {\it Errata},
{\it Ann. of Math.} \textbf{109} (1979), 621--623.

\medbreak\noindent[Takayama2003] Shigeharu Takayama,
Iitaka's fibrations via multiplier ideals. Trans. Amer. Math. Soc. 355 (2003), no. 1, 37--47

\medbreak\noindent[Tsuji2000] Hajime Tsuji, Numerical trivial fibrations, arXiv:math/0001023, 2000.

\medbreak\noindent[Varolin2008]
Dror Varolin,
A Takayama-type extension theorem.
{\it Compos. Math.} \textbf{144} (2008), 522--540.

\medbreak\noindent[Viehweg1980] Eckart Viehweg,
Klassifikationstheorie algebraischer Variet\"aten der Dimension drei.
{\it Compositio Math.} \textbf{41} (1980), 361--400.

\medbreak\noindent[Viehweg-Zuo2003]
Eckart Viehweg and Kang Zuo,
On the Brody hyperbolicity of moduli spaces for canonically polarized manifolds.
{\it Duke Math. J.} \textbf{118} (2003), 103--150.

\medbreak\noindent[Yau1978]
Shing-Tung Yau,
On the Ricci curvature of a compact K\"ahler manifold and the complex Monge-Amp\'ere equation. I.
{\it Comm. Pure Appl. Math.} \textbf{31} (1978), 339--411.

\medbreak\noindent[Weil1958]
Andr\'e Weil,
{\it Introduction \`a l'\'etude des vari\'et\'es k\"ahl\'eriennes}.
Publications de l'Institut de Math\'ematique de l'Universit\'e de Nancago, VI. Actualit\'es Sci. Ind. no.1267, Hermann, Paris 1958.

\bigbreak\noindent{\it Author's mailing address}: Department of
Mathematics, Harvard University, Cambridge, MA 02138, U.S.A.

\medbreak\noindent {\it Author's e-mail address}:
siu@math.harvard.edu

\end{document}